# Different numerical estimators for main effect global sensitivity indices

Sergei Kucherenko[a*], Shufang Song[a]

[a]Centre for Process Systems Engineering, Imperial College London, London, SW7 2AZ, UK

[*]e-mail address: s.kucherenko@imperial.ac.uk

**Abstract.** The variance-based method of global sensitivity indices based on Sobol' sensitivity indices became very popular among practitioners due to its easiness of interpretation. For complex practical problems computation of Sobol' indices generally requires a large number of function evaluations to achieve reasonable convergence. Four different direct formulas for computing Sobol' main effect sensitivity indices are compared on a set of test problems for which there are analytical results. These formulas are based on high-dimensional integrals which are evaluated using MC and QMC techniques. Direct formulas are also compared with a different approach based on the so-called "double loop reordering" formula. It is found that the "double loop reordering" (DLR) approach shows a superior performance among all methods both for models with independent and dependent variables.

## 1 Introduction

Global sensitivity analysis (GSA) complements Uncertainty Quantification in that it offers a comprehensive approach to model analysis by quantifying how the uncertainty in model output is apportioned to the uncertainty in model inputs [1,2]. Unlike local sensitivity analysis, GSA estimates the effect of varying a given input (or set of inputs) while all other inputs are varied as well, thus providing a measure of interactions among variables. GSA is used to identify key parameters whose uncertainty most affects the output. This information then can be used to rank variables, fix unessential variables and decrease problem dimensionality. The variance-based method of global sensitivity indices based on Sobol' sensitivity indices became very popular among practitioners due to its efficiency and easiness of interpretation [3,4]. There are two types of Sobol' sensitivity indices: the main effect indices, which estimate the individual contribution of each input parameter or a group of inputs to the output variance, and the total sensitivity indices, which measure the total contribution of a single input factor or a group of inputs including interactions with all other inputs [5].

For models with independent variables there are efficient direct formulas which allow to compute Sobol' indices directly from function values. These formulas are based on high-dimensional integrals which can be evaluated via MC/QMC techniques [1,4,5]. For complex practical problems computation of Sobol' indices generally requires a large number of function evaluations to achieve reasonable convergence. More efficient formulas for evaluation of Sobol' main effect indices using direct integral formulas were suggested by Saltelli [6]. Kucherenko *et al* [7,8] developed further Saltelli's approach by suggesting new formula which significantly improves the computational accuracy of Sobol' main effect indices with small values. Sobol' and Myshetskaya [9] and Owen [10] suggested their versions of improved direct formulas. In this work we compare original and existing improved direct formulas. For models with dependent inputs we consider a novel approach for estimation Sobol' indices developed in [11]. We also compare direct formulas using MC estimators based on MC and QMC sampling with the so-called double loop approach on a set of test problems which for which there are analytical results for the values of Sobol' indices. The double loop approach has been discarded in the past as being inefficient in comparison with direct formulas but due to the improvements in the algorithms suggested by Plischke [12] it became an interesting alternative to the direct formulas. Further we call this approach as "double loop reordering" (DLR).

Evaluation of Sobol' main effect indices remains to be an active area of research: we could mention application of RBD [13], various metamodelling methods [14,15,16] and some other attempts to improve direct formulas [17]. We also note that a new method for improving the efficiency of the Monte Carlo estimates for the Sobol' total sensitivity indices was developed in [18].

This paper is organized as follows. The next section introduces ANOVA decomposition and Sobol' sensitivity indices. In Section 3 we present different estimators of the main effect sensitivity indices. Comparison of the efficiency of different estimators is considered in Section 4. Finally, conclusions are given in Section 5.

## 2 Sobol' sensitivity indices

Consider the square integrable function $f(x)$ defined in the unit hypercube $H^d = [0,1]^d$. The decomposition of $f(\boldsymbol{x})$

$$f(x) = f_0 + \sum_{i=1}^{d} f_i(x_i) + \sum_{i=1}^{d}\sum_{j>i}^{d} f_{ij}(x_i, x_j) + \cdots + f_{12\cdots d}(x_1, \cdots, x_d), \quad (1)$$

where

$$f_0 = \int_{H^d} f(x)dx$$

is called ANOVA if conditions

$$\int_{H^d} f_{i_1...i_s} dx_{i_k} = 0$$

are satisfied for all different groups of indices $i_1,...,i_s$ such that $1 \le i_1 < i_2 < ... < i_s \le d$. These conditions guarantee that all terms in (1) are mutually orthogonal with respect to integration [4].

The variances of the terms in the ANOVA decomposition add up to the total variance:

$$D = \int_{H^d} f^2(x)dx - f_0^2 = \sum_{s=1}^{d} \sum_{i_1<\cdots<i_s}^{d} D_{i_1...i_s},$$

where components $D_{i_1...i_s} = \int_{H^s} f^2_{i_1...i_s}(x_{i_1},...,x_{i_s})dx_{i_1}...dx_{i_s}$ are called partial variances.

Main effect global sensitivity indices are defined as ratios

$$S_{i_1...i_s} = D_{i_1...i_s} / D.$$

Further we will consider sensitivity indices for a single index:

$$S_i = D_i / D.$$

Total partial variances account for the total influence of the factor $x_i$:

$$D_i^{tot} = \sum_{<i>} D_{i_1...i_s},$$

where the sum $\sum_{<i>}$ is extended over all different groups of indices $i_1,...,i_s$ satisfying condition $1 \le i_1 < i_2 < ... < i_s \le d$, $1 \le s \le d$, where one of the indices is equal $i$ [1,4]. The corresponding total sensitivity index is defined as

$$S_i^{tot} = D_i^{tot} / D.$$

Sobol' also introduced sensitivity indices for subsets of variables [3,4]. Consider two complementary subsets of variables $y$ and $z$:

$$x = (y, z).$$

Let $y = (x_{i_1},...,x_{i_m})$, $1 \le i_1 < ... < i_m \le n$, $K = (i_1,...,i_m)$. Here $m$ is a cardinality of a subset $y$. The variance corresponding to $y$ is defined as

$$D_y = \sum_{s=1}^{m} \sum_{(i_1<\cdots<i_s)\in K} D_{i_1...i_s}. \qquad (2)$$

$D_y$ includes all partial variances $D_{i_1}$, $D_{i_2}$,..., $D_{i_1...i_s}$ such that their subsets of indices

$(i_1,...,i_s) \in K$.

The total variance $D_y^{tot}$ is defined as

$$D_y^{tot} = D - D_z$$

$D_y^{tot}$ consists of all $\sigma_{i_1...i_s}^2$ such that at least one index $i_p \in K$ while the remaining indices can belong to the complementary to $K$ set $\bar{K}$. The corresponding Sobol' sensitivity indices are defined as

$$S_y = D_y / D,$$
$$S_y^{tot} = D_y^{tot} / D.$$

Denote $x_{\sim i} = (x_1,...,x_{i-1},x_{i+1},...,x_d)$ the vector of all variables but $x_i$, then $x \equiv (x_i, x_{\sim i})$ and $f(x) \equiv f(x_i, x_{\sim i})$. The first order component in ANOVA decomposition (1) can be found as

$$f_i(x_i) = \int_{H^d} f(x) dx_{\sim i} - f_0.$$

Then

$$D_i = \int_{H^d} [f_i(x_i)]^2 dx_i = \int_{H^d} \left[ \int_{H^d} f(x) dx_{\sim i} - f_0 \right]^2 dx_i,$$

from which it follows that

$$D_i = \int_{H^d} \left[ \int_{H^d} f(x) dx_{\sim i} \right]^2 dx_i - f_0^2. \tag{3}$$

This formula is used to derive a MC estimator known as the brute force estimator or the double loop method.

There is another approach to derive Sobol' sensitivity indices. If we consider $x$ as a random variable uniformly defined in $H^d$ then $D_i$ can be expressed as [1]:

$$D_i = Var_i[E_{\sim i}(f(x_i, x_{\sim i}) | x_i)]. \tag{4}$$

This representation can be used to derive an extension of Sobol' sensitivity indices for the case of models with dependent variables [11].

## 3 Different formulas and estimators of the main effect sensitivity indices

### 3.1 Original Sobol' formula

Sobol' suggested the following Monte Carlo algorithm for the estimation of $S_y = D_y/D$ [3,4].

Given $x$ and $x'$ being two independent sample points $x = (y, z)$ and $x' = (y', z')$, $D_y$ defined in (2) is calculated using the following formula:

$$D_y = \int f(x) f(y, z') dx dz' - f_0^2. \tag{5}$$

In this case, the Monte Carlo estimator for (5) has a form:

$$D_y \approx \frac{1}{N} \sum_{k=1}^{N} f(y, z) f(y, z') - \left[ \frac{1}{N} \sum_{k=1}^{N} f(y, z) \right]^2, \tag{6}$$

where $N$ is a number of sampled points.

### 3.2 Improved formula of Kucherenko

Kucherenko *et al* [7,8] proposed a new formula for sensitivity indices for sensitivity indices which is especially efficient in the case of indices with small values. Kucherenko and independently Saltelli [6] noticed that $f_0^2$ in (5) can be computed as

$$f_0^2 = \int f(x) f(x') dx dx'. \tag{7}$$

Substituting (7) into (5) and taking out a common multiplier $f(x)$, one can obtaine a new integral representation for $D_y$:

$$D_y = \int f(x) \, f(y, z') - f(x') \, dx dx' \tag{8}$$

and the corresponding Monte Carlo estimator:

$$D_y \approx \frac{1}{N} \sum_{k=1}^{N} f(y, z) \, f(y, z') - f(y', z') \, . \tag{9}$$

Further we refer to this formula as "S-K".

### 3.3 Improved formula of Owen

Owen extended the idea of Kucherenko by using three independent sample points $x = (y, z)$, $x' = (y', z')$ and $x'' = (y'', z'')$ and replacing $f(x)$ by $f(x) - f(y'', z)$ in (8) [10]. As a result $D_y$ is calculated using the following formula:

$$D_y = \int f(x) - f(y'', z) \, f(y, z') - f(x') \, dx dx' dx''. \tag{10}$$

The Monte Carlo algorithm for (10) has a form:

$$D_y \approx \frac{1}{N} \sum_{k=1}^{N} f(y,z) - f(y'',z) \quad f(y,z') - f(y',z') \ . \tag{11}$$

Further we refer to this formula as "Owen".

### 3.4 Improved formula of Sobol' and Myshetzskay

Sobol' and Myshetzskay [9] have argued that formula (8) can be further improved by replacing $f(x)$ in formula by $f(x) - f_0$. Thus, $D_y$ is calculated by using the following formula:

$$D_y = \int f(x) - f_0 \quad f(y,z') - f(x') \, dxdx'dx'' \ . \tag{12}$$

The corresponding Monte Carlo estimator has a form:

$$D_y \approx \frac{1}{N} \sum_{k=1}^{N} f(y,z) - f_0 \quad f(y,z') - f(y',z') \ . \tag{13}$$

Following Owen's classification we further call this formula "Oracle".

### 3.5 Double loop reordering approach

Formula (3) can be used to derive the double loop (the brute force) MC estimator. In this case $N$ points $x^{(j)}$, $j = 1, 2, ..., N$ are generated from the joint probability distribution (PDF) We consider the cases of models with independent and dependent inputs. For each random variable $y = x_i$, the sample set $x^{(j)}$, $j = 1, 2, ..., N$ is sorted and subdivided in $M$ equally populated partitions (bins) each containing $N_m = N/M$ points ($M < N$). Within each bin we calculate the local mean value $E_Z \ f(y,z) | y \approx \frac{1}{N_m} \sum_{k=1}^{N_m} f(y_k, z_k)$. Finally, the variance of all conditional averages is computed as

$$D_y \approx \frac{1}{M} \sum_{j=1}^{M} \left( \frac{1}{N_m} \sum_{k=1}^{N_m} f(y_k^j, z_k^j) \right)^2 - f_0^2 \ . \tag{14}$$

The subdivision in bins is done in the same way for all inputs using the same set of sampled points. This approach we further call the double loop reordering (DLR). A critical issue is the link between $N$ and $M$. It was suggested in [12] to use as a "rule of thumb" $M \sim \sqrt{N}$.

In this work we used Sobol' sequences for QMC sampling [19,20]. To preserve their uniformity properties $N$ should always be equal to $N = 2^p$, where $p$ is an integer. It makes observing the "rule of thumb" more challenging. We used dependence of $M$ and $N_m$ versus $N$ shown in Fig.

1.

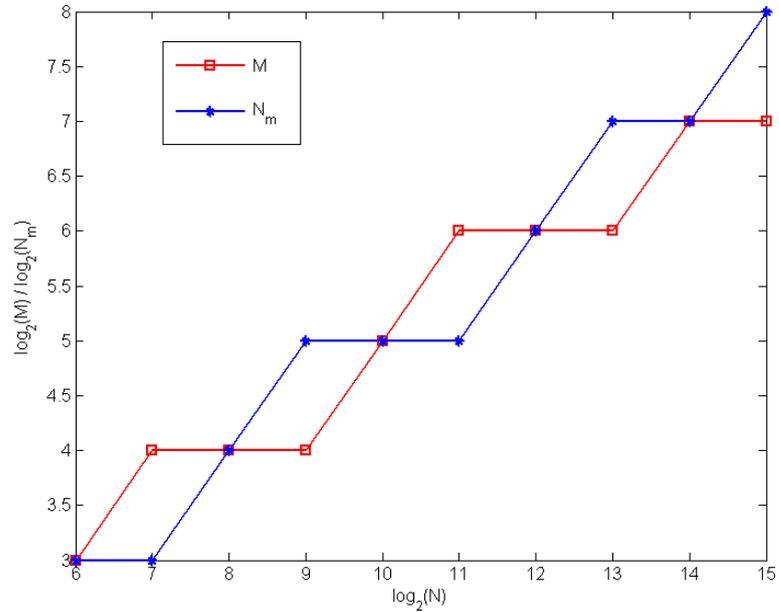

Fig. 1. Dependence of the number of partitions (bins) $M$ and sampled points in each partition $N_m$ versus $N$

We note that although it is possible to extend application of DLR from a single index ($m=1$) to the case of two indices ($m=2$), its extension to $m$ higher than 2 is not practical. Another limitation of DLR in that there is no similar "brute force" formula which allows to compute total Sobol' sensitivity indices.

### 3.6 Number of function evaluations

The five considered MC estimators converge to the same values of the main effect sensitivity indices, but the number of function evaluations per one $i$'th input for each of these methods is different. Table 1 shows the number of function evaluations $N_{CPU}$ required to compute the whole set of sensitivity indices $\{S_i, S_y^{tot}\}$ for a $d$ dimensional function $f(x)$ with independent inputs. Here $N$ is a number of sampled points. We also included $N_{CPU}$ for metamodel based computation of sensitivity indices [15,16]. Here HDMR stands for high dimensional model representation [16].

Table 1: Number of required function evaluations $N_{CPU}$

| Method | Sobol' | S-K | Owen | Oracle | DLR | HDMR |
|---|---|---|---|---|---|---|
| Number of function evaluations $N_{CPU}$ | $N(2d+1)$ | $N(d+2)$ | $N(2d+2)$ | $N(d+2)$ | $N$ | $N$ |

For models with dependent inputs the number of function evaluations $N_{CPU}$ required to compute the whole set of sensitivity indices $\{S_i, S_y^{tot}\}$ $N_{CPU} = N(2d+2)$ [11].

## 4 Numerical tests

The objective of this section is to compare performances of MC and QMC estimators of considered formulas for main effect Sobol' sensitivity indices $S_i$ for models with independent inputs, *i.e.* direct Sobol' formula, Sobol'-Kucherenko (S-K) formula, Owen's formula, Oracle's formula and DLR on a set of test cases for which analytical values of Sobol' sensitivity indices are known. For models with dependent inputs a formula from [11] was also compared with DLR.

For studying the accuracy, the root mean square error (RMSE) $\varepsilon$ is determined using $K$ independent runs:

$$\varepsilon_i(N) = \left[\frac{1}{K}\sum_{k=1}^{K}(S_i^{(n),k} - S_i^{(a)})^2\right]^{1/2}, \quad (15)$$

where $S_i^{(n)}$ and $S_i^{(a)}$ are the numerical and analytical values of $S_i$. Numerical values $S_i^{(n)}$ are computed at $N$, which is reflected in the dependence $\varepsilon_i(N)$. For the MC method all runs are statistically independent. For QMC integration for each run a different part of the Sobol' sequence was used. For all tests we took $K=10$.

The QMC convergence rate is known to be $\varepsilon_{QMC} = \dfrac{O(\ln N)^d}{N}$ [20]. In pracitce, the rate of convergence for QMC methods appears to be approximately equal to $O(N^{-\alpha})$, with $0.5 \leq \alpha \leq 1$. For the MC method $\alpha = 0.5$. QMC method in most cases outperforms MC in terms of convergence [8]. In practical tests the RMSE $\varepsilon_i(N)$ is approximated by the formula $cN^{-\alpha}$, $0 < \alpha < 1$, and the convergence rate $\alpha$ is extracted from fitted trend lines. We consider convergence rates for various estimators versus $N$ and $N_{CPU}$.

### 4.1 Models with independent inputs

**Test 1:** Consider a model

$$f(x) = a_1 x_1 + a_2 x_2 + \ldots + a_d x_d,$$

where $x_1, x_2, \ldots, x_d$ are independent normal variables: $x_i \sim N(\mu_i, \sigma_i^2)$, $a_1, a_2, \ldots, a_d$ are constant coefficients. The PDF of the output $Y$ is normally distributed, *i.e.*

$$Y \sim N\left(\sum_{i=1}^{d} a_i \mu_i, \sum_{i=1}^{d} a_i^2 \sigma_i^2\right),$$

while the PDF of the conditional output is

$$Y \mid X_i \sim N\left(a_i x_i + \sum_{j=1, j\neq i}^{d} a_j \mu_j, \sum_{j=1, j\neq i}^{d} a_j^2 \sigma_j^2\right).$$

It's easy to see that the analytical values of Sobol' indices are $S_i = S_i^{tot} = \dfrac{a_i^2 \sigma_i^2}{\sum_{j=1}^{d} a_j^2 \sigma_j^2}$.

We consider the case $d=4$ with the mean values and standard deviations $\mu = (1, 3, 5, 7)$ and $\sigma = (1, 1.5, 2, 2.5)$, respectively with all coefficients $a_i = 1$ ($i=1, 2, 3, 4$). The analytical values of $S_i$ are $S_i = \{0.0741, 0.167, 0.296, 0.463\}$. It is clear from the convergence and RMSE plots presented in Figs. 2-4 that

1) For the MC method for input $i = 1$ which has a small value $S_1$, all three improved formulas have much higher convergence rate than the original Sobol' formula with Owen's formula outperforming all other methods (Fig. 3). DLR has a similar performance to Oracle's formula. Situation is different for input $i = 4$ which has rather high value $S_4$: although all three improved formulas have higher convergence rate than the original Sobol' formula but the difference between the original Sobol' and S-K formulas are smaller. Owen and Oracle formulas are the most efficient among all direct formulas but the clear winner is DLR. We also note, that the extracted convergence rate $\alpha$ is close to 0.5 as expected for the MC method.

2) For the QMC method for input $i = 1$ (a small value $S_1$), all three improved formulas have much higher convergence rate than the original Sobol' formula with Owen and Oracle formulas outperforming all other methods (Fig. 2, 4). DLR shows the highest performance superior to direct formulas. For input $i = 4$ three improved formulas have higher convergence rate than the original Sobol' formula but the difference between the original Sobol' formula and S-K and Owen's formulas are smaller similarly to the previous case with MC sampling. Oracle's formula is the most efficient among all direct formulas but DLR exhibits even higher convergence. The extracted convergence rate $\alpha$ is close to 1.0 as expected for the QMC method in the case of low effective dimension [8].

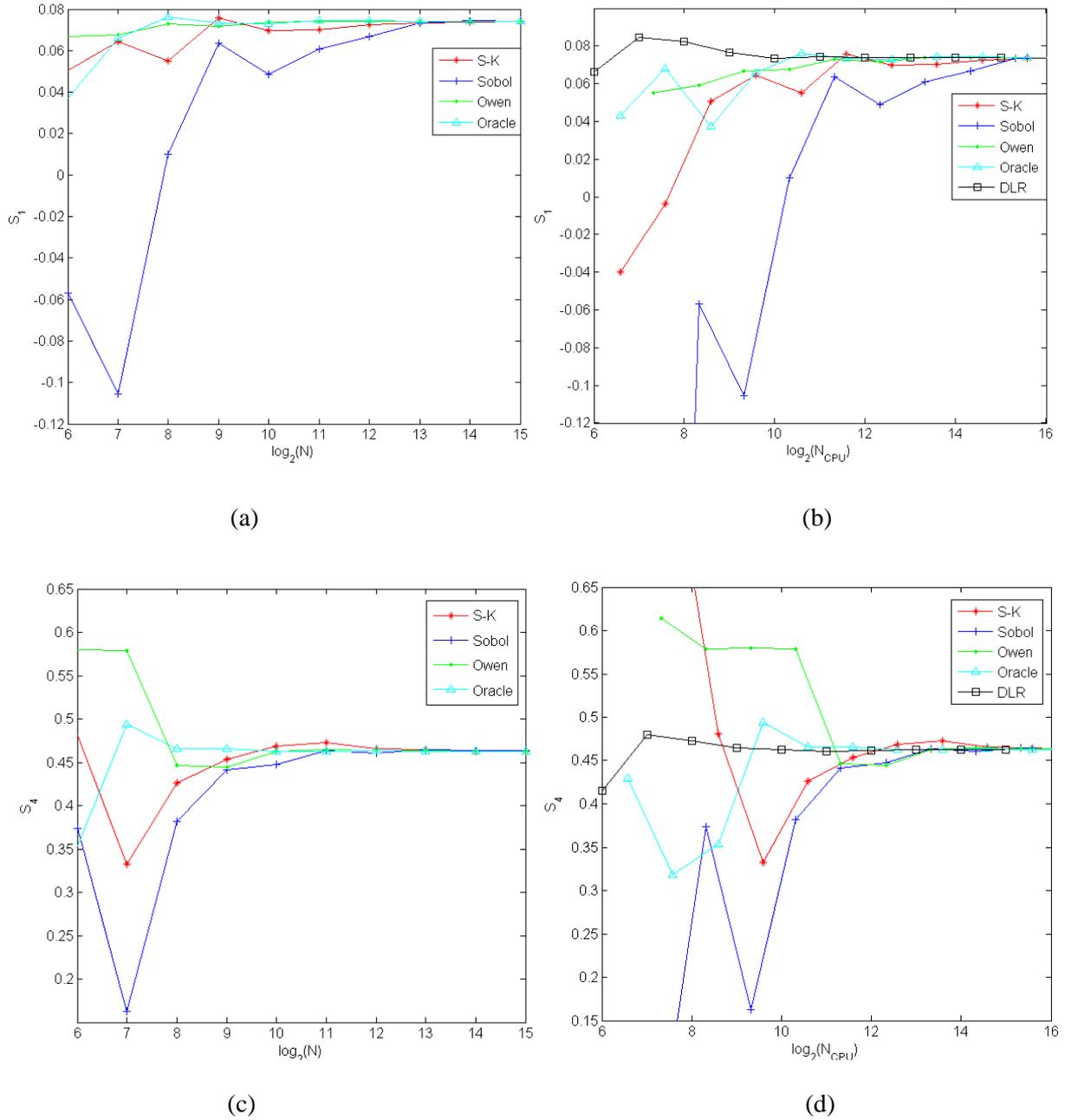

Fig. 2. Test case 1: Convergence plots of $S_i$, $i$ = 1 (a), (b), $i$ = 4 (c), (d). QMC sampling. The red line refers to S-K formula; the blue line refers to Sobol' formula, the green line refers to Owen's formula, the cyan line refers to Oracle formula, the black line refers to DLR. On the left: (a), (c) the values of $S_i$ obtained at the same number of $N$. On the right: (b), (d) the values of $S_i$ obtained at the same number of $N_{CPU}$.

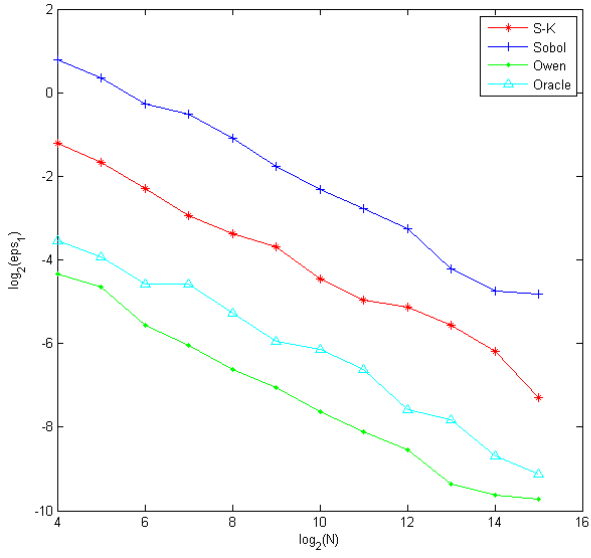
(a)

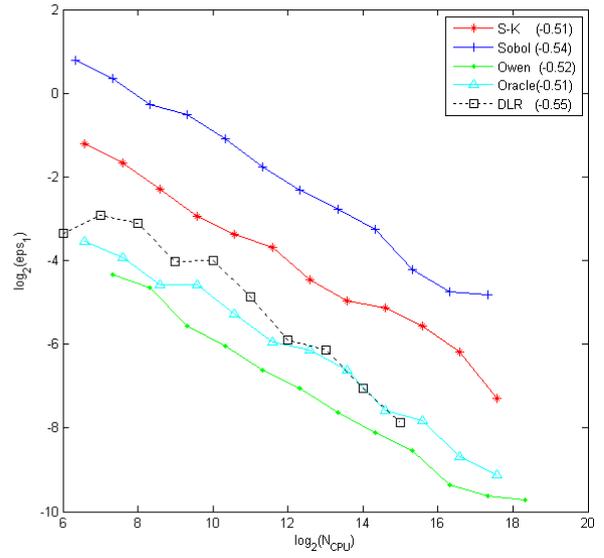
(b)

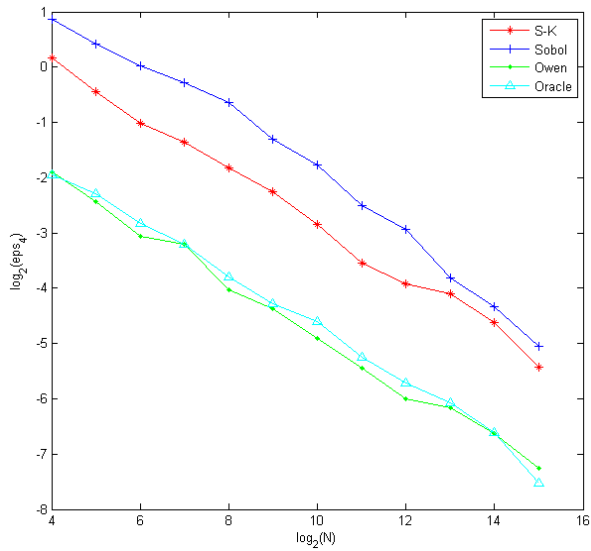
(c)

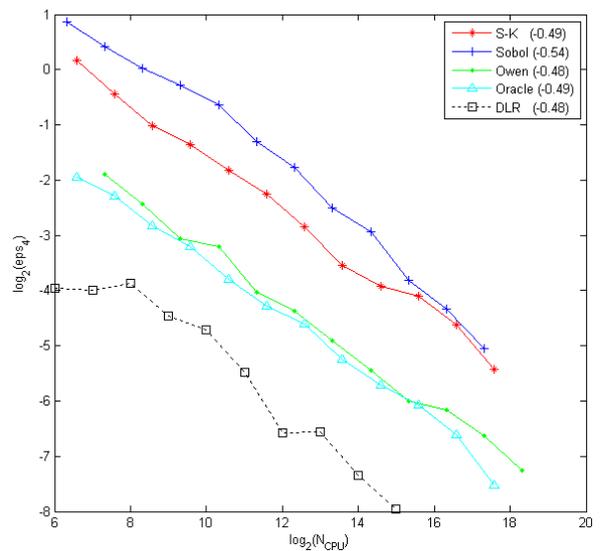
(d)

Fig. 3. Test case 1: The RMSE $\varepsilon_i$ $i = 1$ (a), (b), $i = 4$ (c), (d) versus the number of $N$ (on the left: (a), (c)) and the number of $N_{CPU}$ (on the right: (b), (d)). MC sampling. The red line refers to S-K formula; the blue line refers to Sobol' formula, the green line refers to Owen's formula, the cyan line refers to Oracle formula, the black line refers to DLR.

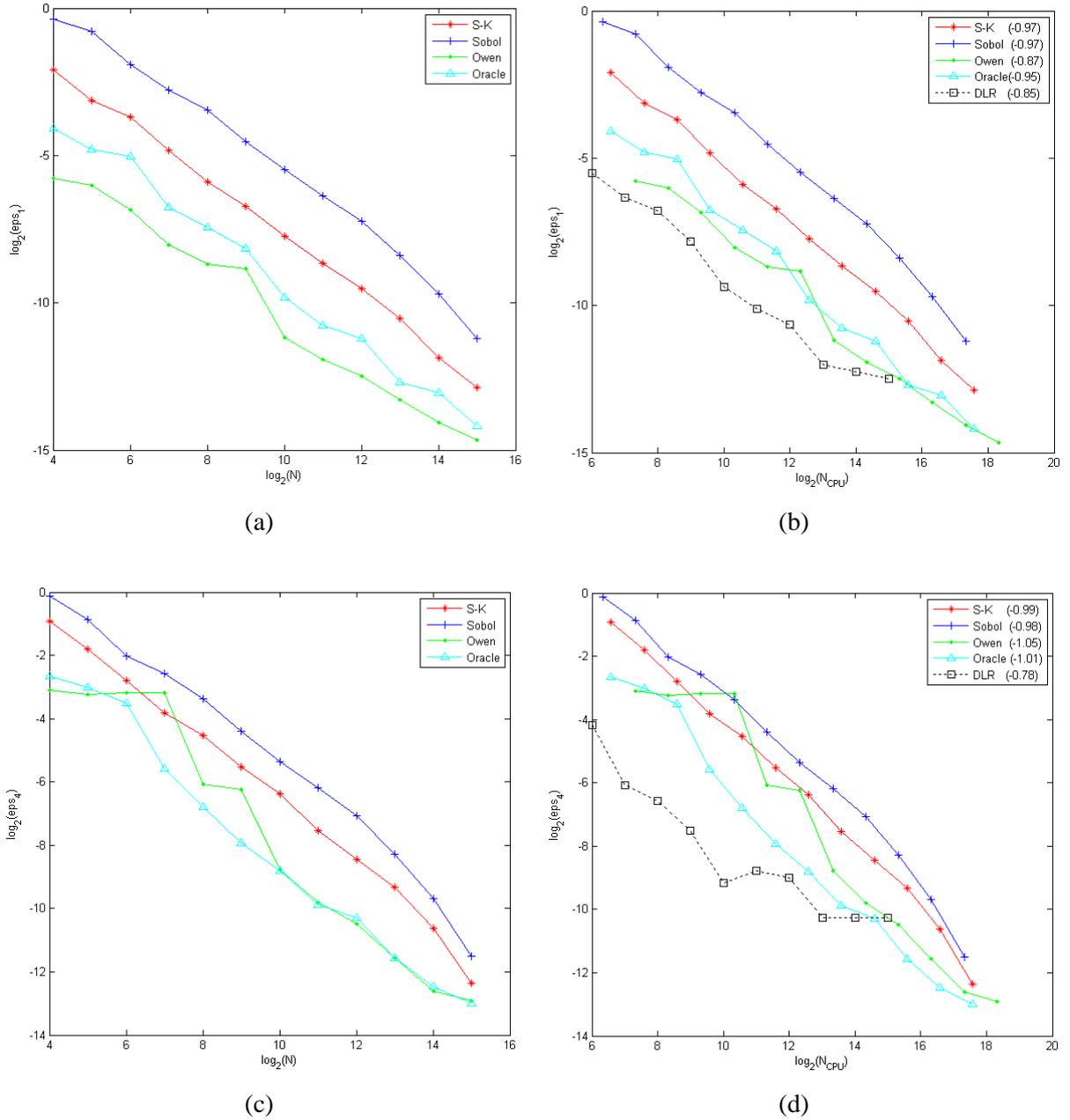

Fig. 4. Test case 1: The RMSE $eps_i$ $i = 1$ (a), (b), $i = 4$ (c), (d) versus the number of $N$ (on the left: (a), (c)) and the number of $N_{CPU}$ (on the right: (b), (d)). QMC sampling. The red line refers to S-K formula; the blue line refers to Sobol' formula, the green line refers to Owen's formula, the cyan line refers to Oracle formula, the black line refers to DLR.

**Test 2:** Consider a model

$$f(x) = x_1x_3x_5 + x_1x_3x_6 + x_1x_4x_5 + x_1x_4x_6 + x_2x_3x_4 + x_2x_3x_5 + x_2x_4x_5 + x_2x_5x_6 + x_2x_4x_7 + x_2x_6x_7$$

in which all seven variables are independent lognormal with the mean values 2, 3, 0.001, 0.002, 0.004,

0.005 and 0.003, respectively. All the standard deviations are all equal to 0.4214. This model was taken from [21]. Values of $S_i$ are $S_i$ = {0.0350, 0.330, 0.0157, 0.0857, 0.174, 0.221, 0.0477}. From the convergence and RMSE plots presented in Figs. 5-7 we can conclude that

1) For the MC method for input $i = 2$ which has a large value $S_2$, all three improved formulas have a higher convergence rate than the original Sobol' formula with Owen and Oracle formulas outperforming all other methods (Fig. 6). DLR has a superior performance over other methods. Situation is different for input $i = 4$ which has small $S_4$: all three improved formulas have a much higher convergence rate than the original Sobol' formula. Owen and Oracle formulas are the most efficient among all direct formulas and they show a similar performance to DLR.

2) For the QMC method for inputs $i = 2, 6$ (large values $S_i$), all three improved formulas show slightly higher convergence rate than the original Sobol' formula (Fig. 5, 7). For input $i = 4$ (small $S_4$) three improved formulas have much higher convergence rate than the original Sobol' formula. DLR shows a superior performance over other methods for all inputs.

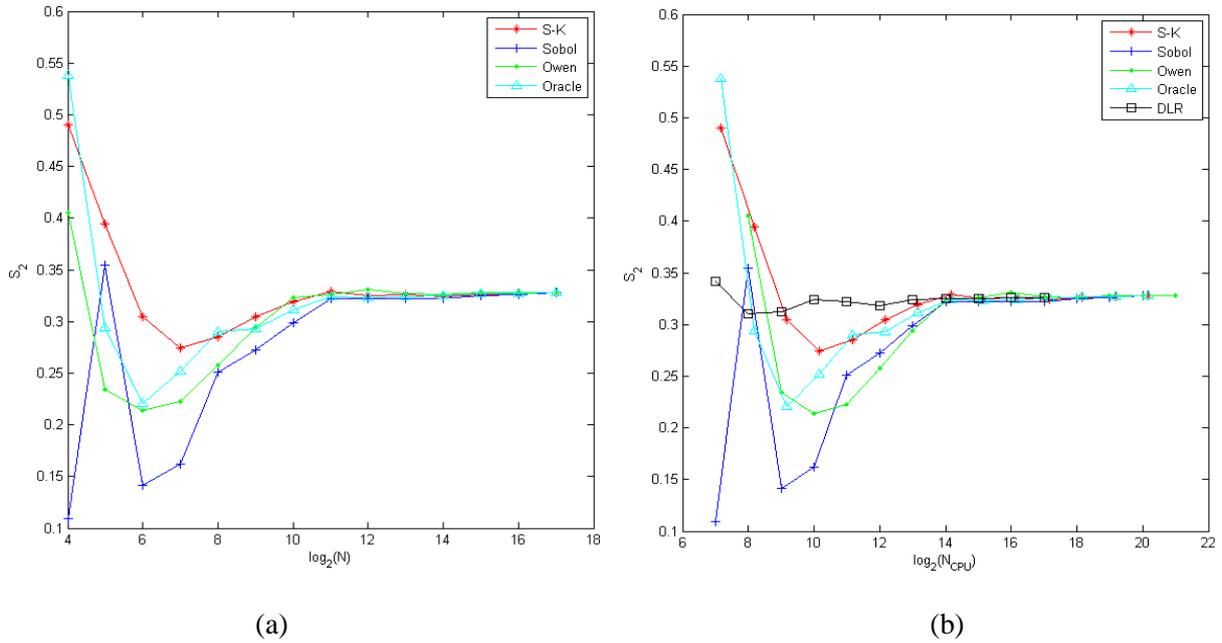

(a)              (b)

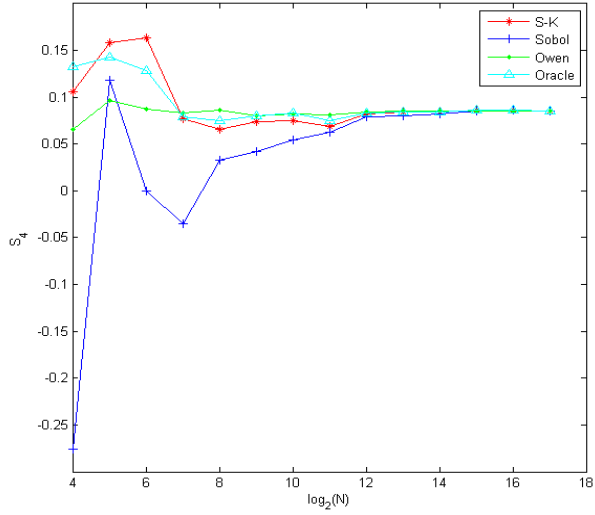
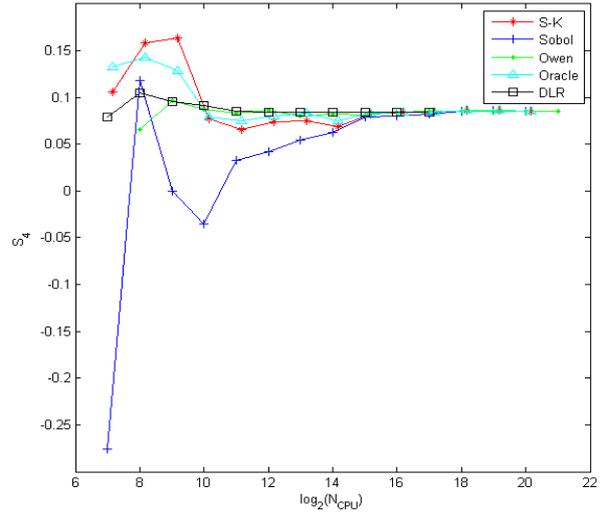

(c) (d)

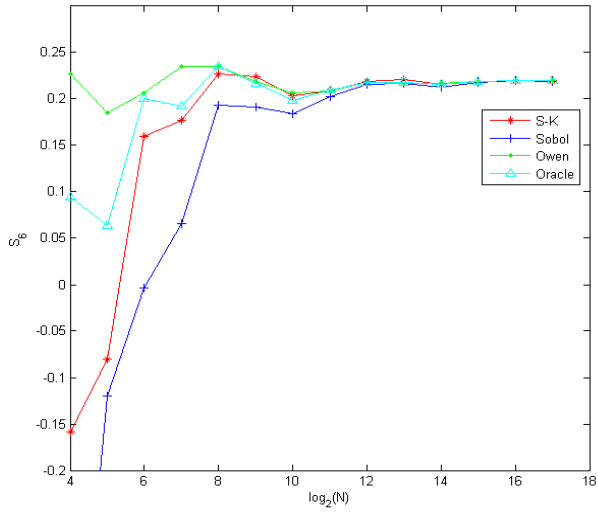
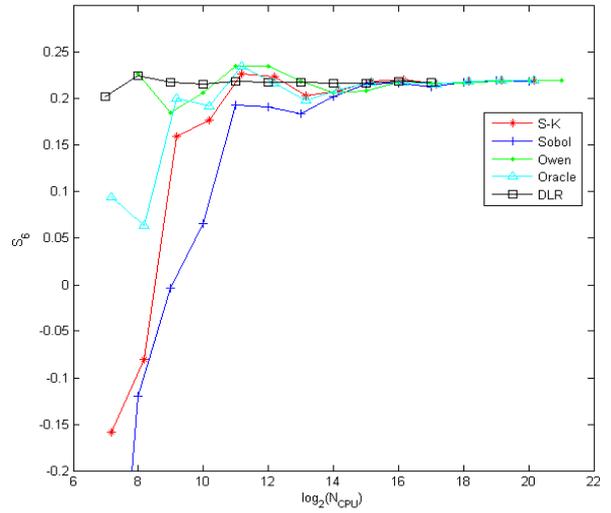

(e) (f)

Fig. 5. Test case 2: Convergence plots of $S_i$, $i = 2$ (a), (b), $i = 4$ (c), (d), $i = 6$ (e), (f). QMC sampling. The red line refers to S-K formula; the blue line refers to Sobol' formula, the green line refers to Owen's formula, the cyan line refers to Oracle formula, the black line refers to DLR. On the left: (a), (c), (e) the values of $S_i$ obtained at the same number of $N$. On the right: (b), (d), (f) the values of $S_i$ obtained at the same number of $N_{CPU}$.

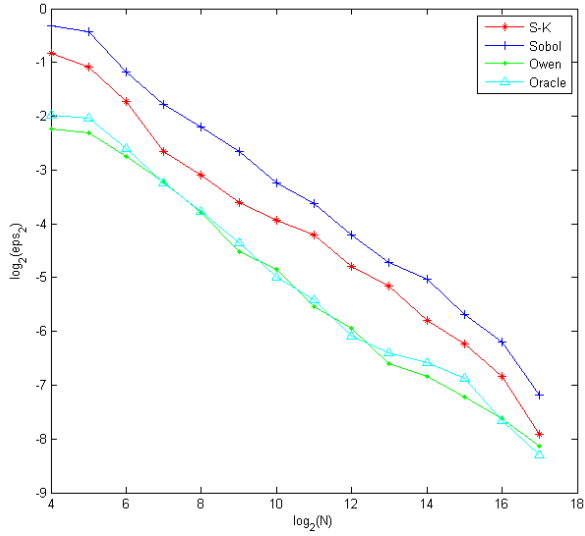

(a)

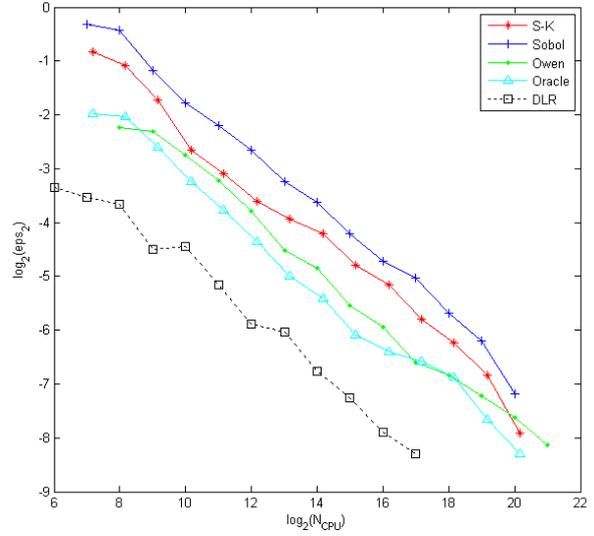

(b)

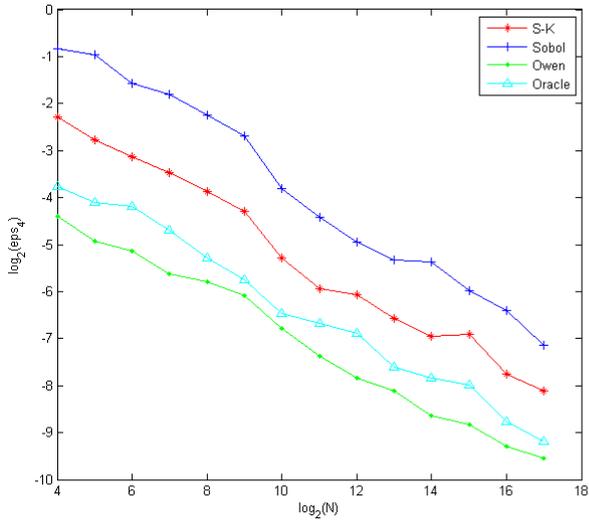

(b)

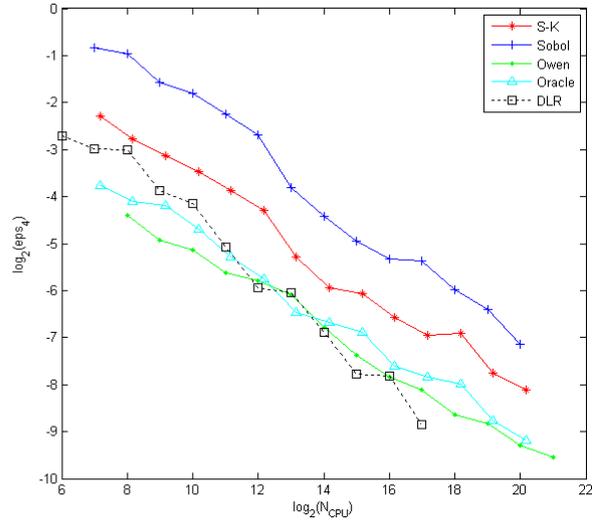

(d)

Fig. 6. Test case 2: The RMSE $\varepsilon_i$ $i = 2$ (a), (b), $i = 4$ (c), (d) versus the number of $N$ (on the left: (a), (c)) and the number of $N_{CPU}$ (on the right: (b), (d)). MC sampling. The red line refers to S-K formula; the blue line refers to Sobol' formula, the green line refers to Owen's formula, the cyan line refers to Oracle formula, the black line refers to DLR.

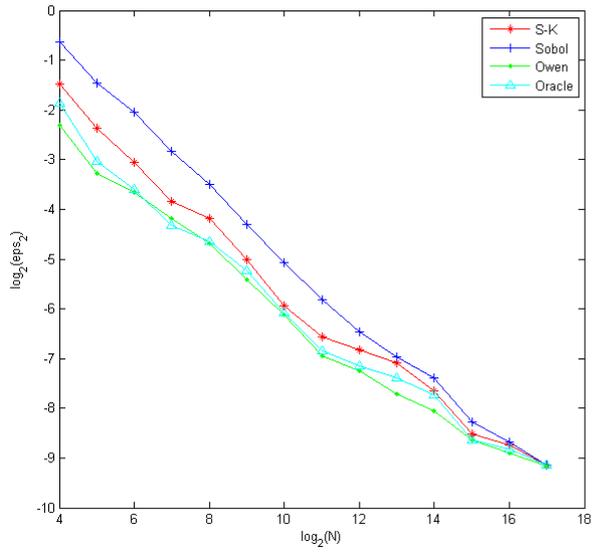 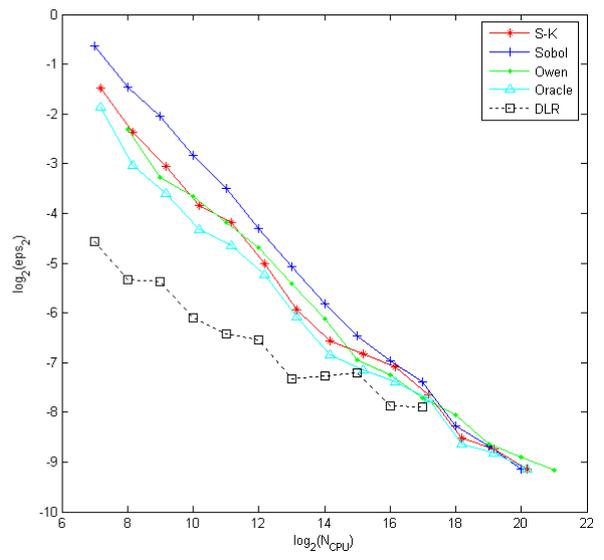

(a)           (b)

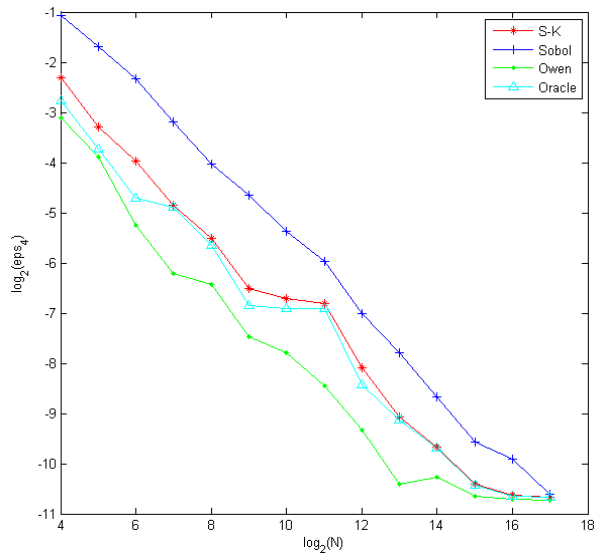 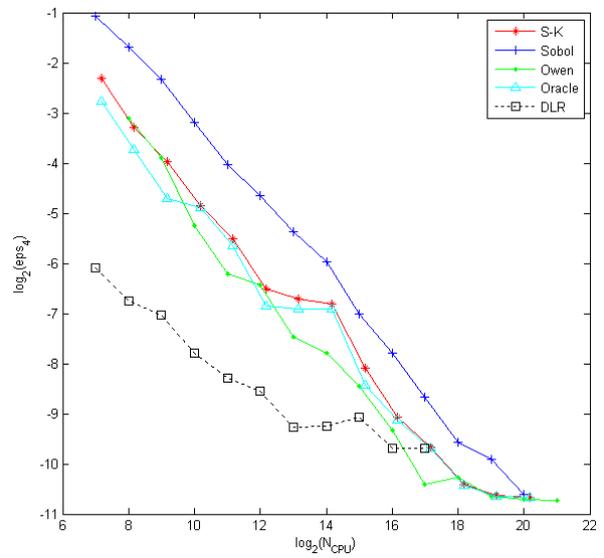

(c)           (d)

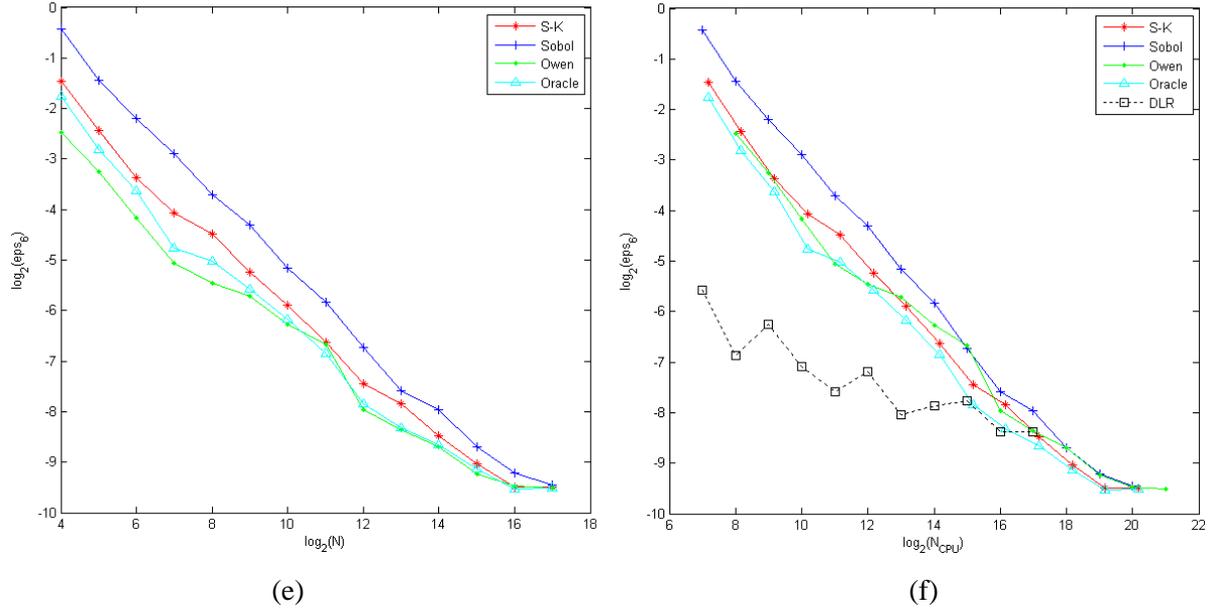

(e)           (f)

Fig. 7. Test case 2: The RMSE $\varepsilon_i$ $i = 2$ (a), (b), $i = 4$ (c), (d), $i = 6$ (e), (f) versus the number of $N$ (on the left: (a), (c), (e)) and the number of $N_{CPU}$ (on the right: (b), (d), (f)). QMC sampling. The red line refers to S-K formula; the blue line refers to Sobol' formula, the green line refers to Owen's formula, the cyan line refers to Oracle formula, the black line refers to DLR.

**Test 3:** The Ishigami function: $f(\boldsymbol{x}) = \sin(x_1) + 7(\sin x_2)^2 + 0.1 x_3^4 \sin(x_1)$ is often used as a benchmark for sensitivity analyses studies [1]. Here $x_i$, $i = 1, 2, 3$ are uniformly distributed on the interval $[-\pi, \pi]$. The Sobol's sensitivity indices have the following values: $S_i$ = {0.314, 0.442, 0.00}. From the convergence and RMSE plots presented in Figs. 8-10 we can conclude that

1) For the MC method for inputs $i = 1, 2$ which have large values of $S_i$ Oracle's formula slightly outperforms other methods (Fig. 9). For input $i = 3$ for which $S_3 = 0.0$ all three improved formulas have a much higher convergence rate than the original Sobol' formula. DLR has a superior performance over other methods for all three inputs.
2) For the QMC method for inputs $i = 1, 2$ direct formulas show a similar convergence rate (Fig. 8, 10). DLR shows higher performance only for inputs $i = 1$. Figs. 10, d also presents the results for $S_2$ obtained using the QMC-HDMR method from [16]. Clearly, the QMC- HDMR method shows somewhat better converence than other methods.

.

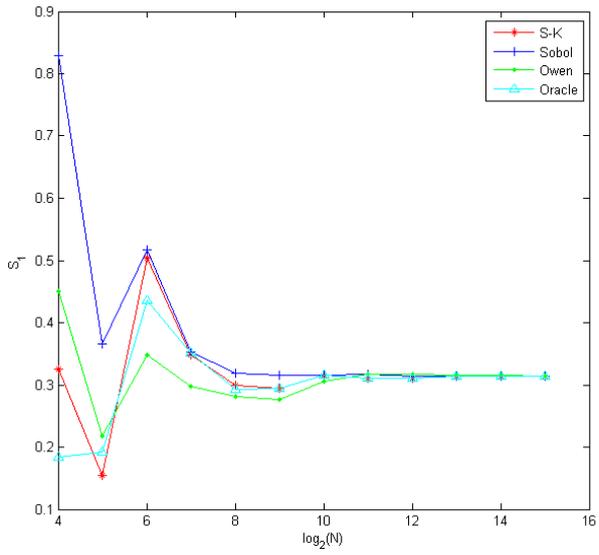 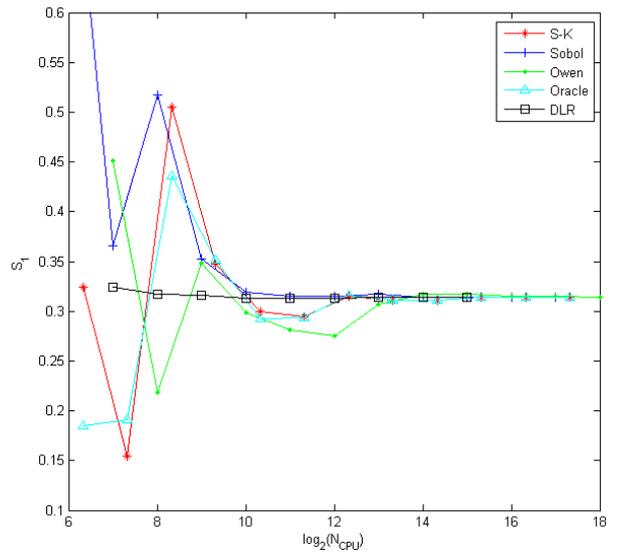

(a) (b)

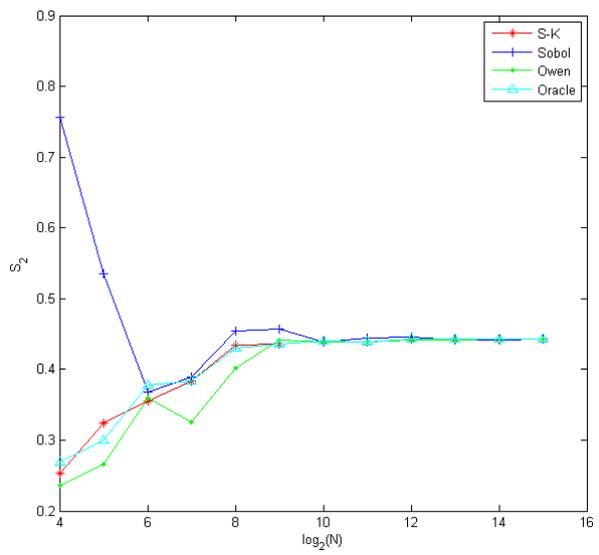 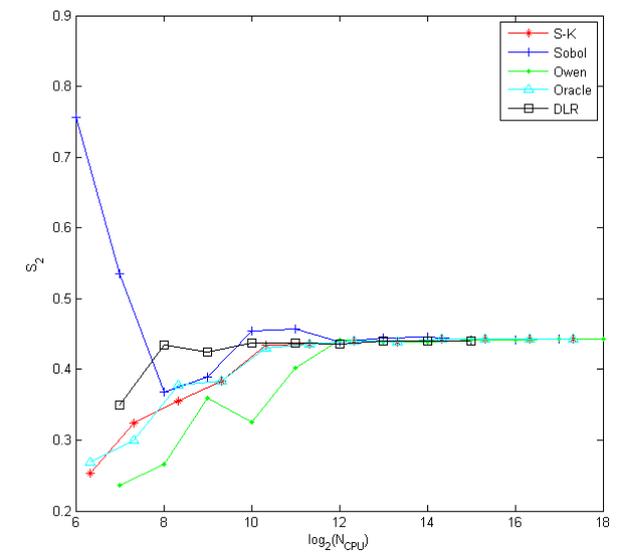

(c) (d)

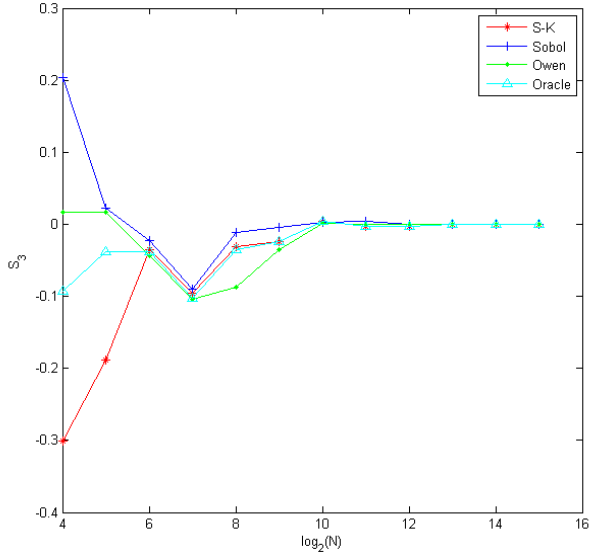
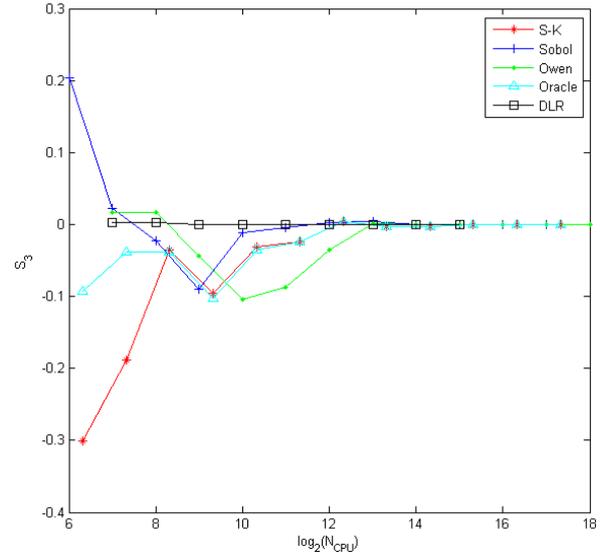

(e)                                               (f)

Fig. 8. Test case 3: Convergence plots of $S_i$, $i = 1$ (a), (b), $i = 4$ (c), (d), $i = 6$ (e), (f). QMC sampling. The red line refers to S-K formula; the blue line refers to Sobol' formula, the green line refers to Owen's formula, the cyan line refers to Oracle formula, the black line refers to DLR. On the left: (a), (c), (e) the values of $S_i$ obtained at the same number of $N$. On the right: (b), (d), (f) the values of $S_i$ obtained at the same number of $N_{CPU}$.

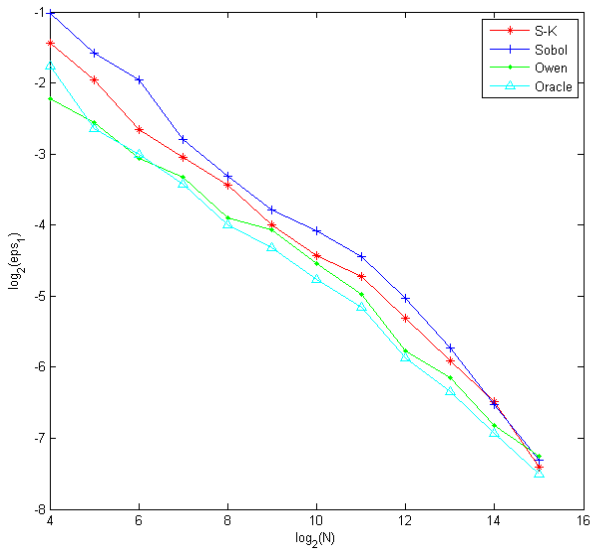
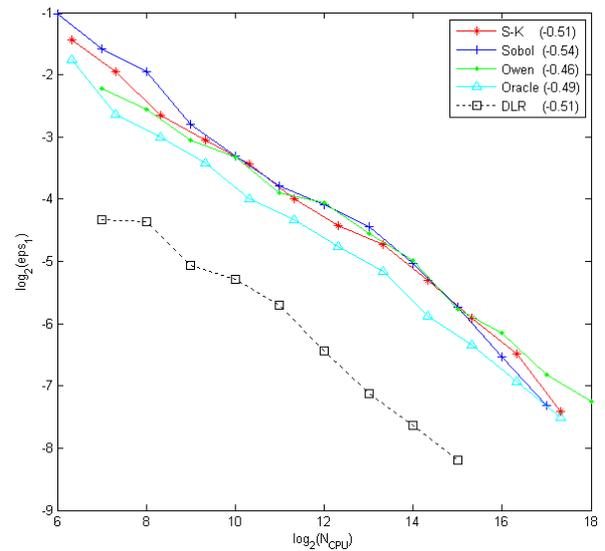

(a)                                               (b)

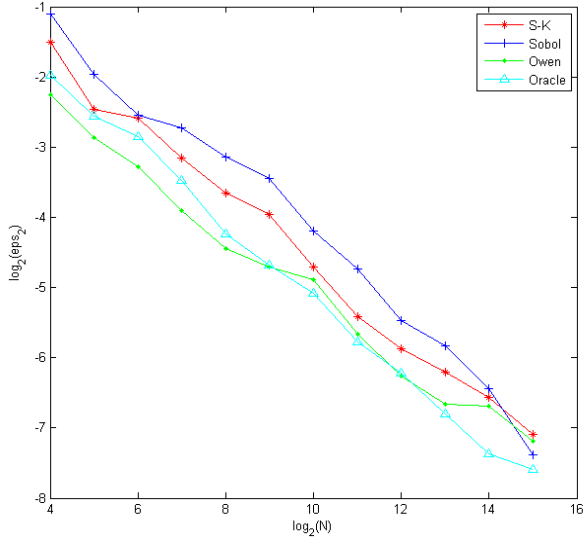
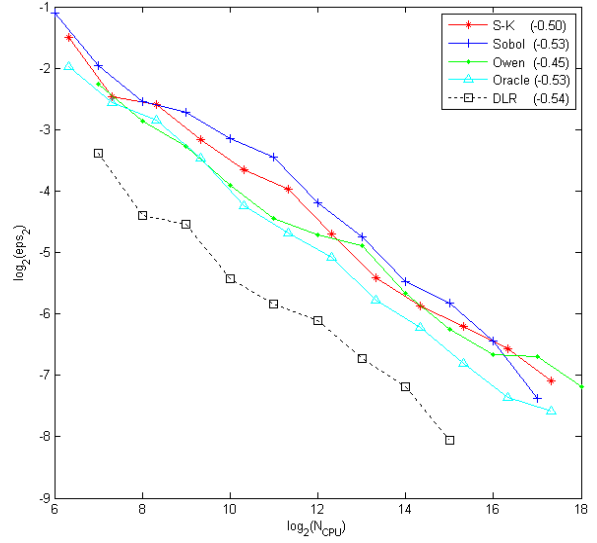

(c)          (d)

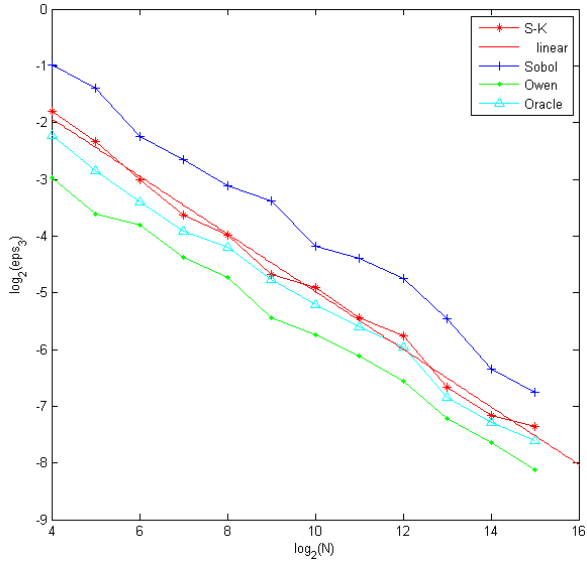
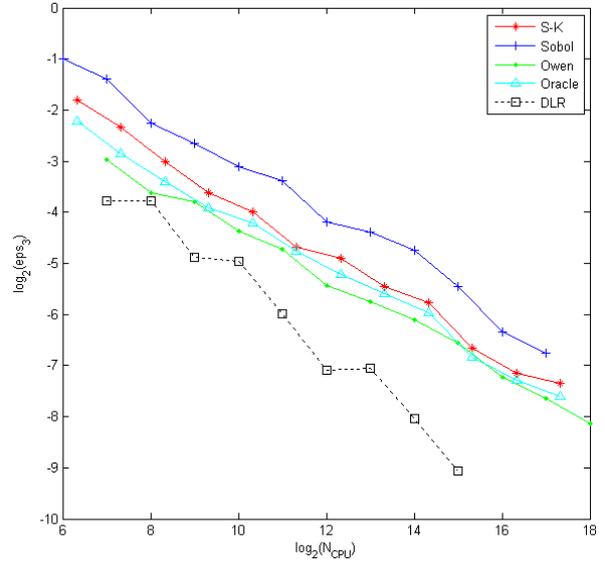

(e)          (f)

Fig. 9. Test case 3: The RMSE $\varepsilon_i$ $i = 1$ (a), (b), $i = 2$ (c), (d), $i = 3$ (e), (f) versus the number of $N$ (on the left: (a), (c), (e)) and the number of $N_{CPU}$ (on the right: (b), (d), (f)). MC sampling. The red line refers to S-K formula; the blue line refers to Sobol' formula, the green line refers to Owen's formula, the cyan line refers to Oracle formula, the black line refers to DLR.

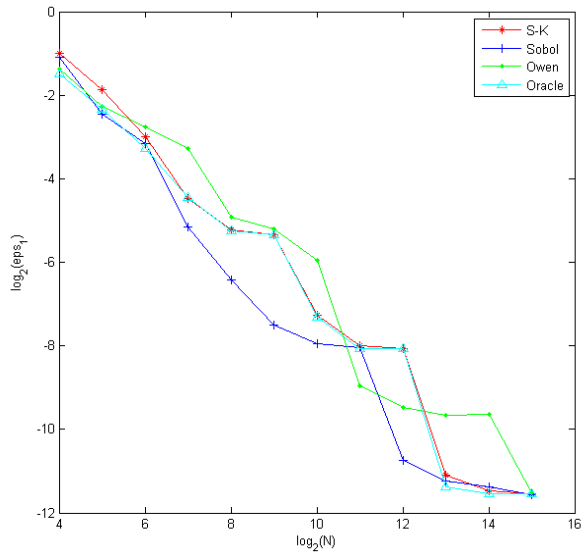

(a)

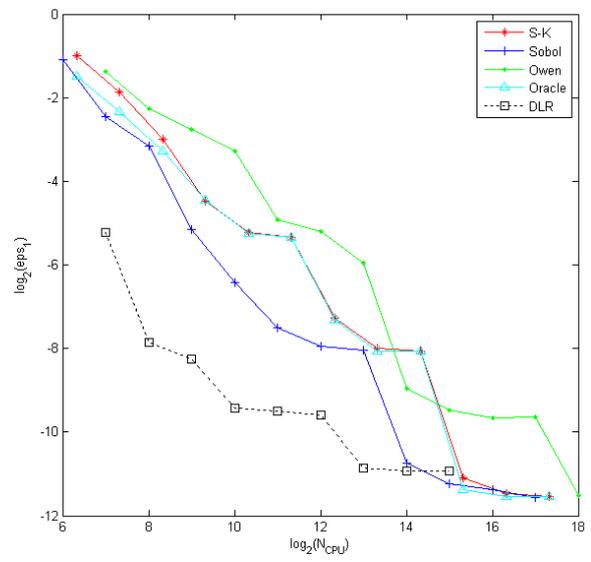

(b)

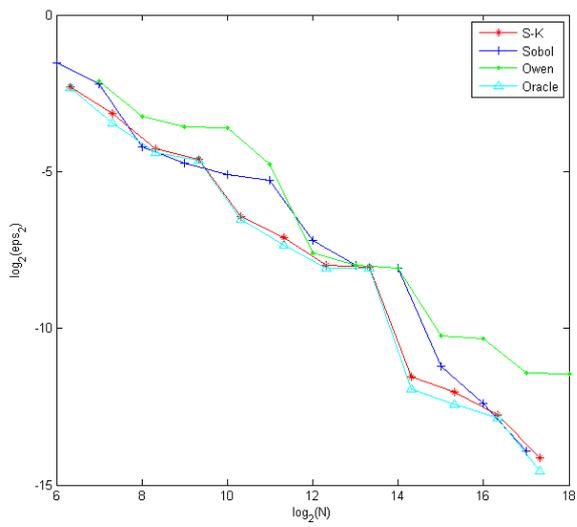

(c)

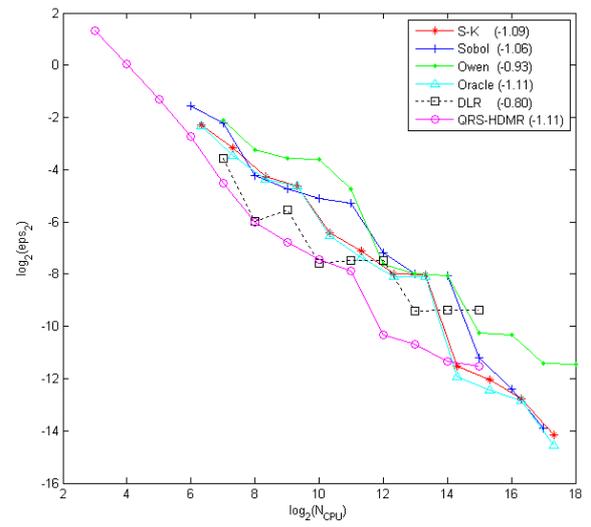

(d)

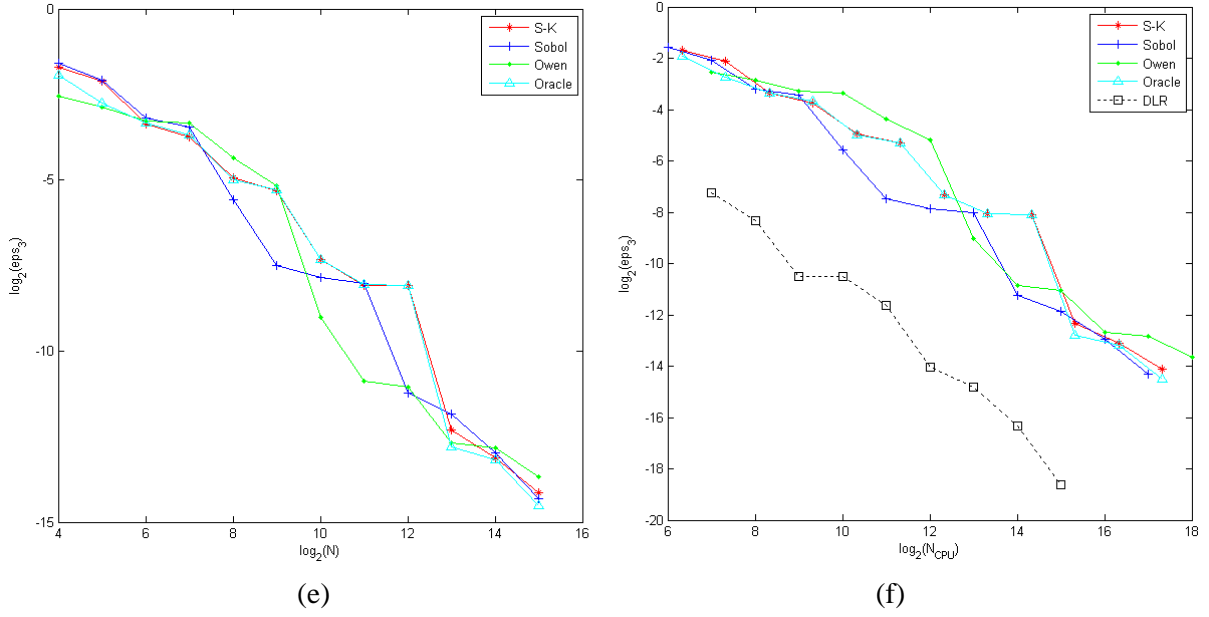

(e)                  (f)

Fig. 10. Test case 3: The RMSE $\varepsilon_i$ $i = 1$ (a), (b), $i = 2$ (c), (d), $i = 3$ (e), (f) versus the number of $N$ (on the left: (a), (c), (e)) and the number of $N_{CPU}$ (on the right: (b), (d), (f)). QMC sampling. The red line refers to S-K formula; the blue line refers to Sobol' formula, the green line refers to Owen's formula, the cyan line refers to Oracle formula, the black line refers to DLR, the magenta line refers to metamodel based computation of $S_i$.

**Test 4:** The g-function

$$f(x) = \prod_{i=1}^{d} g_i = \prod_{i=1}^{d} \frac{|4x_i - 2| + a_i}{1 + a_i},$$

is also often used as a benchmark [4]. Here $d$ is the number of independent input factors ($0 \leq x_i \leq 1$). Parameter $a_i$ determines the importance of the input factor $x_i$.

**Test 4.1:** We consider the 10-dimensional g-function with parameters $a_1 = a_2 = 0$, $a_3 = \cdots = a_{10} = 3$. The analytical values of $S_i$: $S_1 = S_2 = 0.304$, $S_3 = \cdots = S_{10} = 0.019$. For $a_i = 0$ variable $x_i$ is important, for $a_i = 3$ variable $x_i$ is unimportant, hence only the first two variables are important. From the convergence and RMSE plots presented in Figs. 11-13 we can conclude that

1) For the MC method for input $i = 1$ which has a large value of $S_i$, Oracle's formula slightly outperforms other methods (Fig. 12). DLR has a superior performance over other methods. For input $i =$ 3 which has a small value of $S_i$ all three improved formulas have a much higher convergence rate than the original Sobol' formula with Owen's formula outperforming all other methods. DLR shows performance similar to Owen's formula when $N_{CPU} > 2^{11}$.

2) The results for the QMC method qualitatively are similar to those of the MC method (Figs. 11, 13) with the only difference in that DLR shows a superior performance among all other methods for input $i = 3$. However, quantitatively the rate of convergence for the QMC method is much higher than that for the MC method.

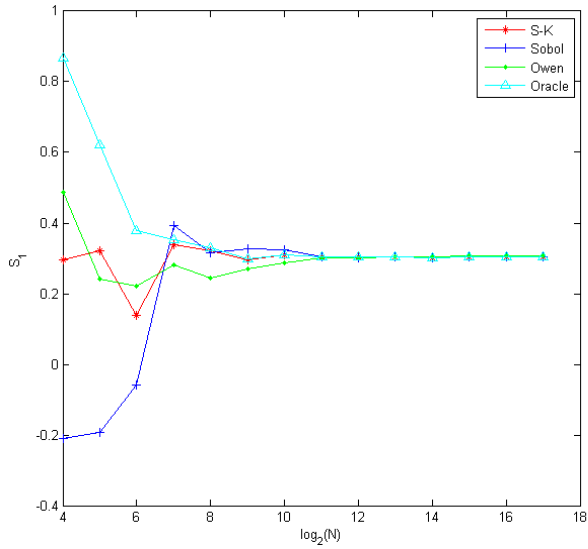
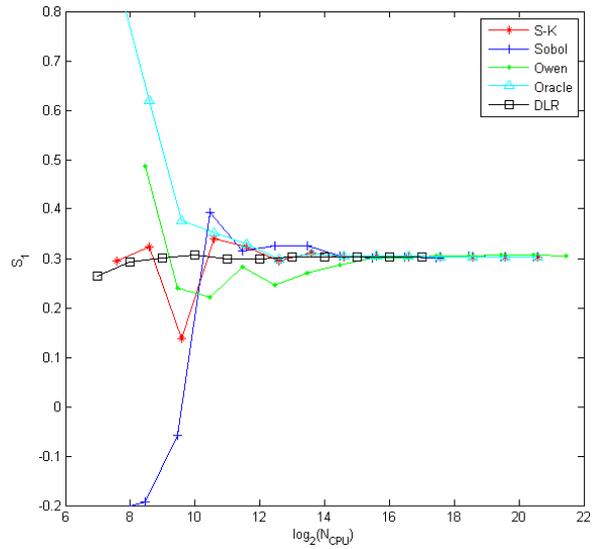

(a) (b)

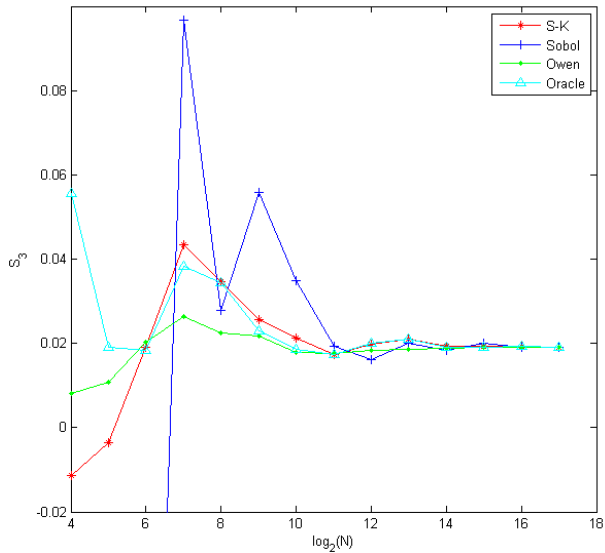
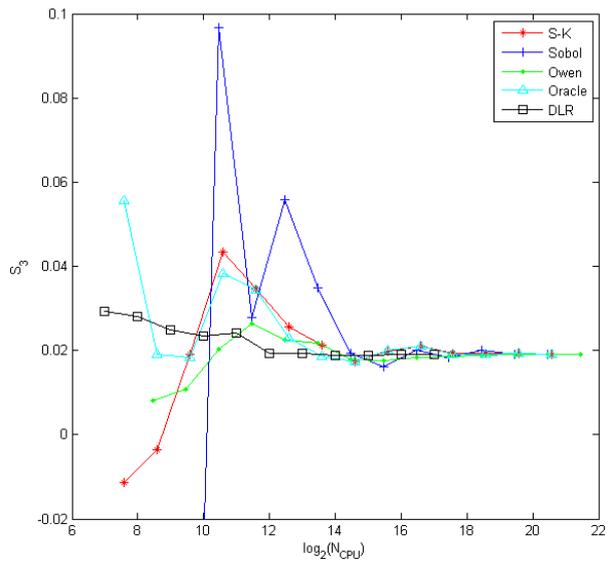

(c) (d)

Fig. 11. Test case 3: Convergence plots of $S_i$, $i = 1$ (a), (b), $i = 3$ (c), (d). QMC sampling. The red line refers to S-K formula; the blue line refers to Sobol' formula, the green line refers to Owen's formula, the cyan line refers to Oracle formula, the black line refers to DLR. On the left: (a), (c) the values of $S_i$ obtained at the same number of $N$. On the right: (b), (d) the values of $S_i$ obtained at the same number of $N_{CPU}$.

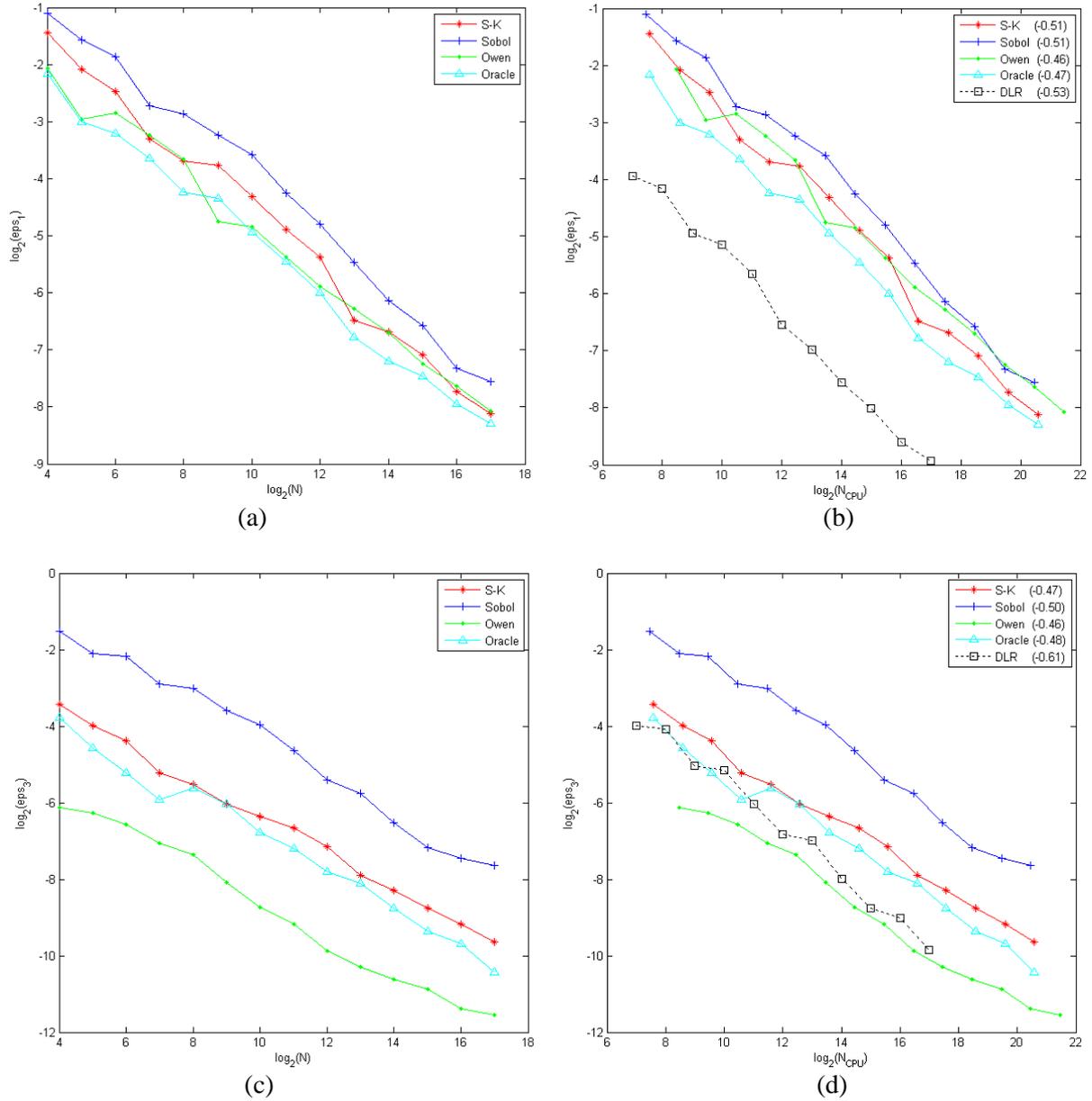

(a)

(b)

(c)

(d)

Fig. 12. Test case 4.1: The RMSE $eps_i$ $i = 1$ (a), (b), $i = 3$ (c), (d) versus the number of $N$ (on the left: (a), (c)) and the number of $N_{CPU}$ (on the right: (b), (d)). MC sampling. The red line refers to S-K formula; the blue line refers to Sobol' formula, the green line refers to Owen's formula, the cyan line refers to Oracle formula, the black line refers to DLR.

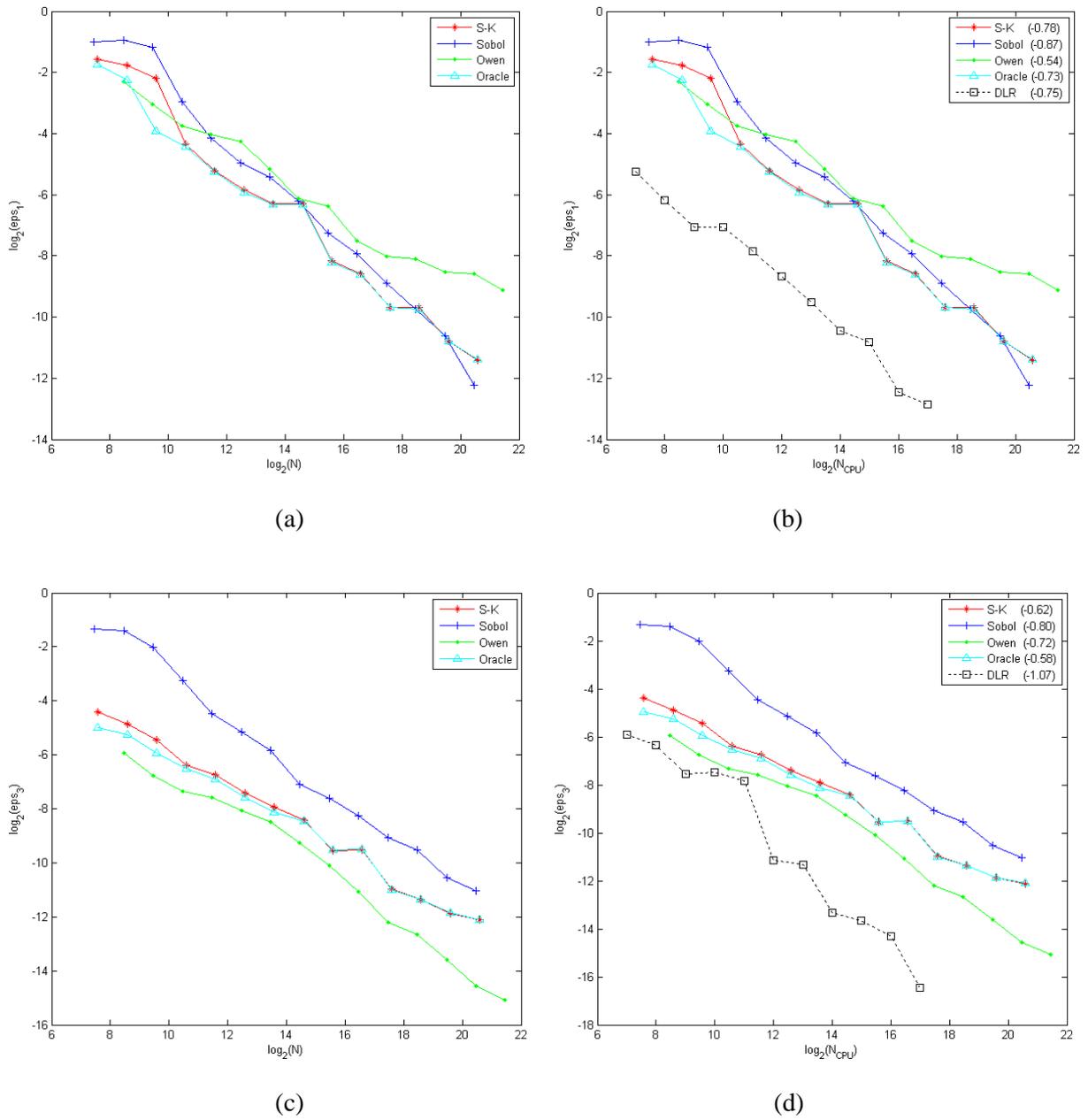

(a)  (b)

(c)  (d)

Fig. 13. Test case 4.1: The RMSE $\varepsilon_i$ $i = 1$ (a), (b), $i = 3$ (c), (d) versus the number of $N$ (on the left: (a), (c)) and the number of $N_{CPU}$ (on the right: (b), (d)). QMC sampling. The red line refers to S-K formula; the blue line refers to Sobol' formula, the green line refers to Owen's formula, the cyan line refers to Oracle formula, the black line refers to DLR.

**Test 4.2**: All parameters $a_i = 0$, $i = 1, 2, ..., 10$, the analytical value of $S_i$ is 0.0199. All inputs are equally important and there are strong interactions between inputs. This is type C function [8] for which the QMC method is not more efficient than MC. From the convergence and RMSE plots presented in Figs. 14-16 we can conclude that

1) For the MC method, all three improved formulas show slightly higher convergence rate than the original Sobol' formula (Fig. 15). DLR outperforms other methods.
2) The results for the QMC method are similar to those of the MC method (Figs. 14, 16) with the only difference in that DLR has a higher higher convergence rate than that for the MC method.

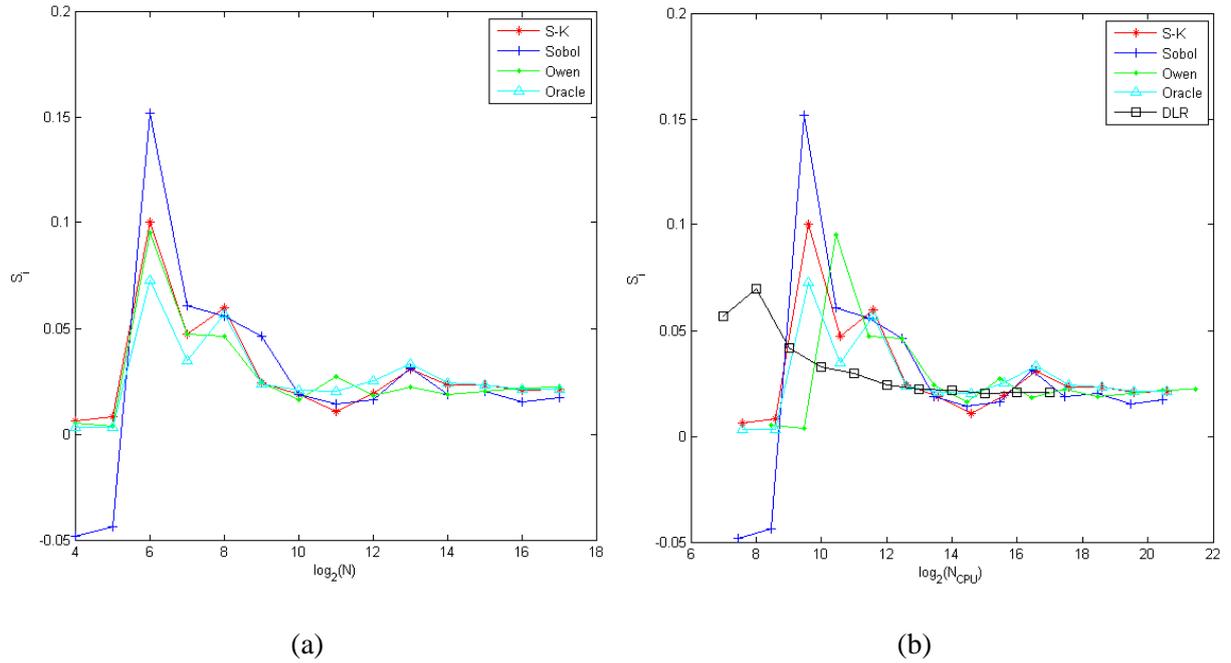

(a)                                      (b)

Fig. 14. Test case 4.2: Convergence plots of $S_i$, $i = 1$. QMC sampling. The red line refers to S-K formula; the blue line refers to Sobol' formula, the green line refers to Owen's formula, the cyan line refers to Oracle formula, the black line refers to DLR. On the left: (a) the values of $S_i$ obtained at the same number of $N$. On the right: (b) the values of $S_i$ obtained at the same number of $N_{CPU}$.

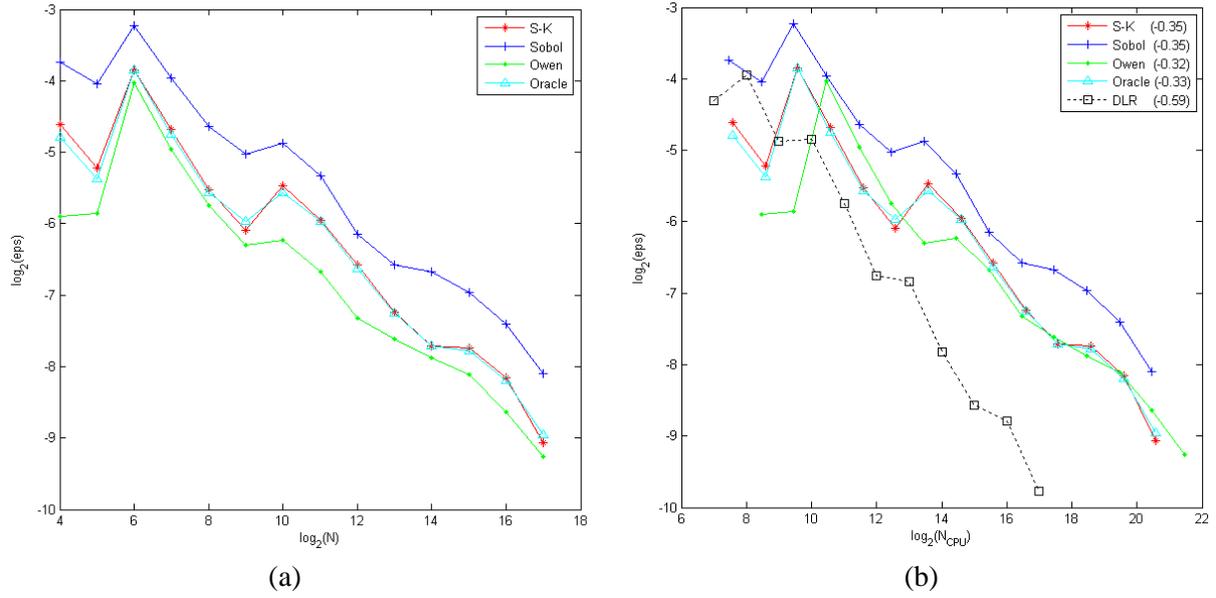

Fig. 15. Test case 4.2: The RMSE $eps_i$ versus the number of $N$ (on the left: (a)) and the number of $N_{CPU}$ (on the right: (b)). MC sampling. The red line refers to S-K formula; the blue line refers to Sobol' formula, the green line refers to Owen's formula, the cyan line refers to Oracle formula, the black line refers to DLR.

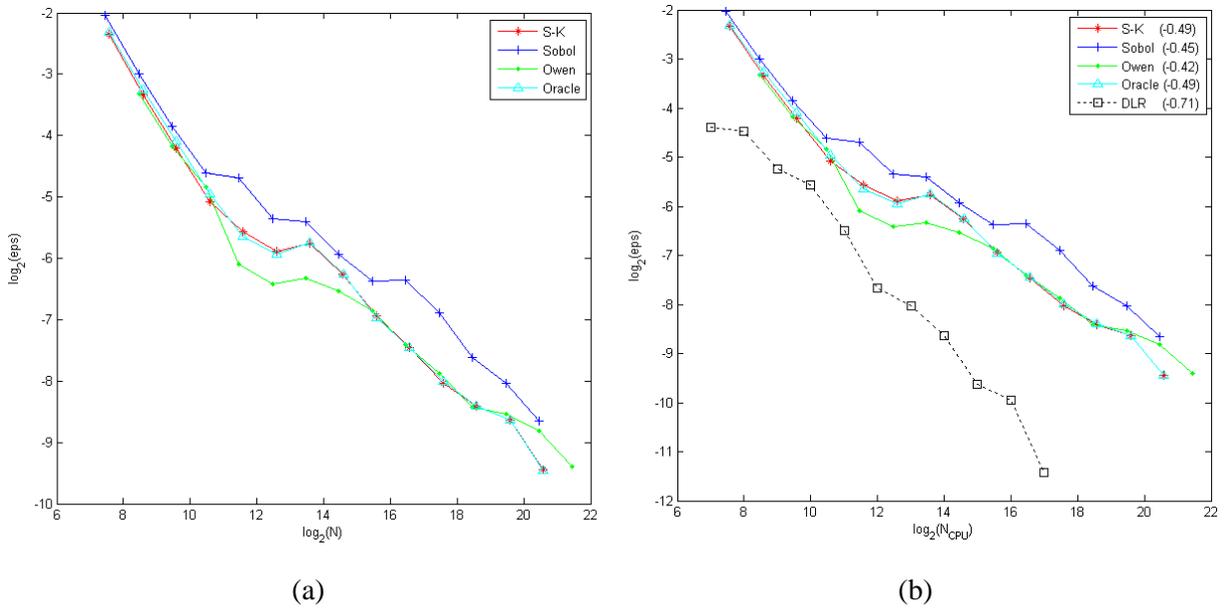

Fig. 16. Test case 4.2: The RMSE $\varepsilon_i$ versus the number of $N$ (on the left: (a)) and the number of $N_{CPU}$ (on the right: (b)). QMC sampling. The red line refers to S-K formula; the blue line refers to Sobol' formula, the green line refers to Owen's formula, the cyan line refers to Oracle formula, the black line refers to DLR.

## 4.2 Models with dependent inputs

**Test case 5:** Consider a model $f(x) = x_1 x_3 + x_2 x_4$, where $(x_1, x_2, x_3, x_4) \sim N(\boldsymbol{\mu}, C_x)$ with $\boldsymbol{\mu} = (0, 0, \mu_3, \mu_4)$ and the covariance matrix

$$C_x = \begin{bmatrix} \sigma_1^2 & \sigma_{12} & 0 & 0 \\ \sigma_{12} & \sigma_2^2 & 0 & 0 \\ 0 & 0 & \sigma_3^2 & \sigma_{34} \\ 0 & 0 & \sigma_{34} & \sigma_4^2 \end{bmatrix}.$$

This test case was considered in [11] were the analytical values of the main and total order indices were presented (Table 2).

Table 2: Test case 5. Analytical values of the main and total order indices.

|  | $x_1$ | $x_2$ | $x_3$ | $x_4$ |
|---|---|---|---|---|
| $S_i$ | $\dfrac{\sigma_1^2 \left(\mu_3 + \mu_4 \rho_{12} \dfrac{\sigma_2}{\sigma_1}\right)^2}{D}$ | $\dfrac{\sigma_2^2 \left(\mu_4 + \mu_3 \rho_{12} \dfrac{\sigma_1}{\sigma_2}\right)^2}{D}$ | 0 | 0 |
| $S_i^{tot}$ | $\dfrac{\sigma_1^2 \left(1 - \rho_{12}^2\right)\left(\sigma_3^2 + \mu_3^2\right)}{D}$ | $\dfrac{\sigma_2^2 \left(1 - \rho_{12}^2\right)\left(\sigma_4^2 + \mu_4^2\right)}{D}$ | $\dfrac{\sigma_1^2 \sigma_3^2 \left(1 - \rho_{34}^2\right)}{D}$ | $\dfrac{\sigma_2^2 \sigma_4^2 \left(1 - \rho_{34}^2\right)}{D}$ |

Here $\rho_{ij} = \dfrac{\sigma_{ij}}{\sigma_i \sigma_j}$ and $D = \sigma_1^2 \left(\sigma_3^2 + \mu_3^2\right) + \sigma_2^2 \left(\sigma_4^2 + \mu_4^2\right) + 2\sigma_{12}(\sigma_{34} + \mu_3 \mu_4)$. For numerical test we used the following parameters:

$$\boldsymbol{\mu} = (0, 0, 250, 400), \quad C_x = \begin{bmatrix} 16 & 2.4 & 0 & 0 \\ 2.4 & 4 & 0 & 0 \\ 0 & 0 & 4 \cdot 10^4 & -1.8 \cdot 10^4 \\ 0 & 0 & -1.8 \cdot 10^4 & 9 \cdot 10^4 \end{bmatrix}.$$

Numerical values of Sobol' sensitivity indices are $S_i = \{0.507, 0.399, 0, 0\}$. From the convergence plots presented in Fig. 17 we can conclude that DLR outperforms the extended version of Sobol' formula [11] for high values of $S_i$ ($i=1, 2$), however for zero values of $S_i$ ($i=3, 4$) it is slightly less efficient than the extended version of Sobol' formula. We note that DLR is much easier to implement algorithmically as it does not require rather complex procedure of sampling from conditional distribution. The CPU time required for the extended version of Sobol' formula is 19.6 s versus only 1.78 s required for DLR.

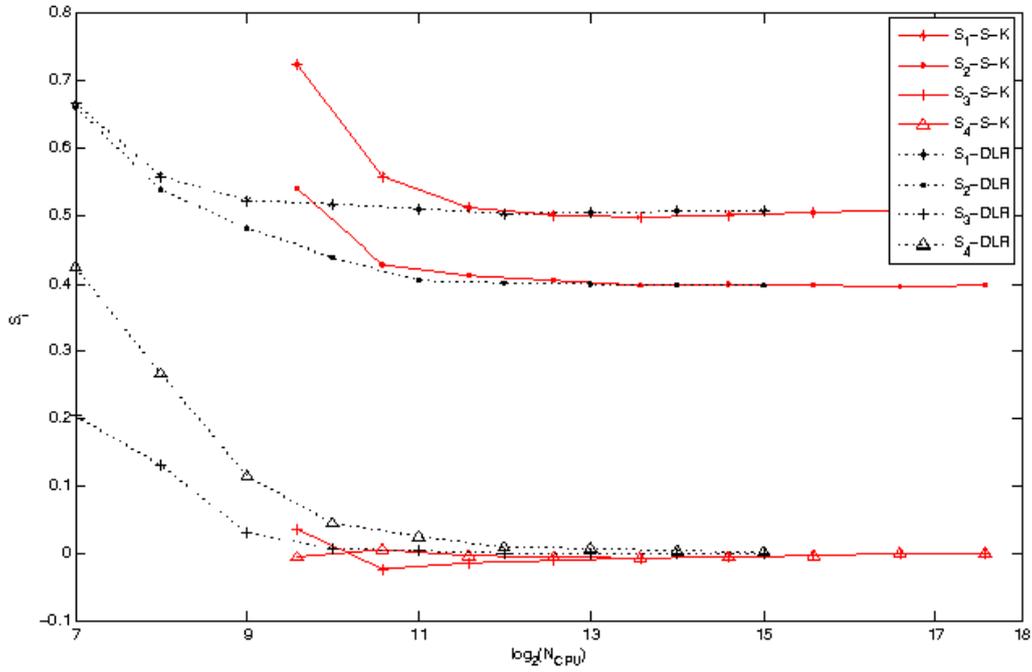

Fig. 17. Test case 5: Convergence plots of $S_i$, $i = 1,..,4$. QMC sampling. The red line refers to S-K formula; the black line refers to DLR. The values of $S_i$ obtained at the same number of $N_{CPU}$. $S_1$ - upper lines with crosses, $S_2$ - lines with circles, $S_3$ - lines with triangles, $S_4$ - lower lines with crosses.

**Test case 6:** Consider the linear model $f(x) = x_1 + x_2 + x_3$, where all input variables are normally distributed with zero mean and the covariance matrix $C_x$:

$$C_x = \begin{bmatrix} 1 & 0 & 0 \\ 0 & 1 & \rho\sigma \\ 0 & \rho\sigma & \sigma^2 \end{bmatrix}.$$

Analytical values of both main effect and total Sobol' sensitivity indices were given in [11] (Table 3).

Table 3: Test case 6. Analytical values of the main and total order indices.

|  | $x_1$ | $x_2$ | $x_3$ |
|---|---|---|---|
| $S_i$ | $\dfrac{1}{2+\sigma^2+2\rho\sigma}$ | $\dfrac{\left(1+\rho\sigma\right)^2}{2+\sigma^2+2\rho\sigma}$ | $\dfrac{\left(\rho+\sigma\right)^2}{2+\sigma^2+2\rho\sigma}$ |
| $S_i^{tot}$ | $\dfrac{1}{2+\sigma^2+2\rho\sigma}$ | $\dfrac{1-\rho^2}{2+\sigma^2+2\rho\sigma}$ | $\dfrac{\left(1-\rho^2\right)\sigma^2}{2+\sigma^2+2\rho\sigma}$ |

For numerical test we used the following parameters: $\boldsymbol{\mu}=(0,\ 0,\ 0), \sigma=2, \rho=-0.8$. Similarly to the previous test case from the convergence plots presented in Fig. 18 we can conclude that DLR outperforms the extended version of Sobol' formula for high values of $S_i$ ($i$=1, 2), however for small values of $S_i$ ($i$=3) both methods show a similar performance.

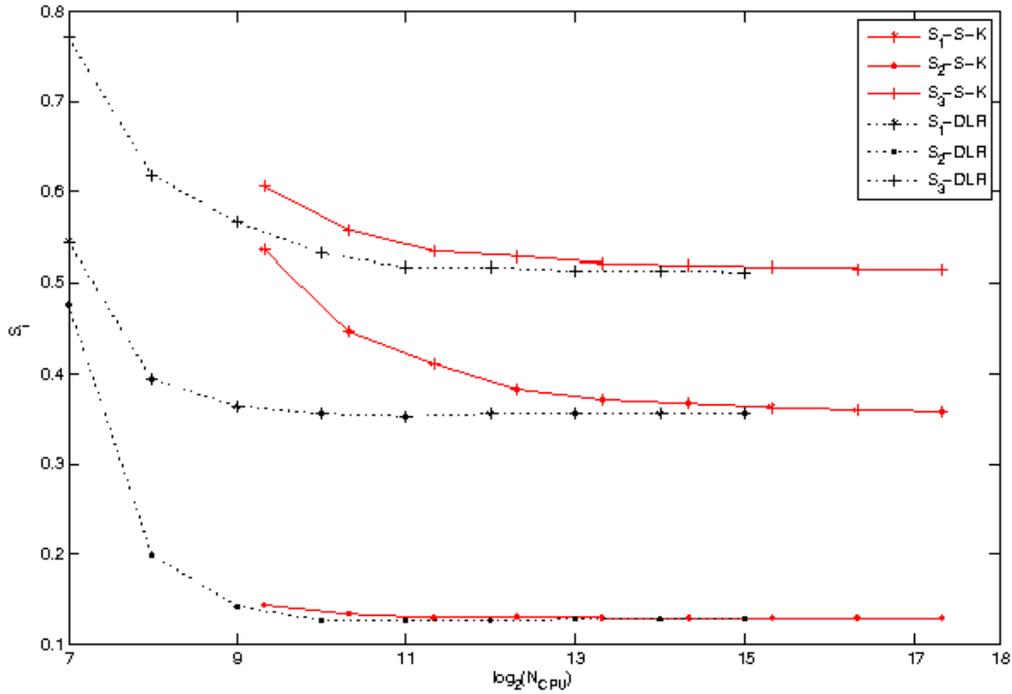

Fig. 18. Test case 6: Convergence plots of $S_i$, $i=1,2,3$. QMC sampling. The red line refers to S-K formula; the black line refers to DLR. The values of $S_i$ obtained at the same number of $N_{CPU}$. $S_1$ - upper lines with crosses, $S_2$ - lines with circles, $S_3$ - medium lines with crosses.

## 5  Conclusions

In this paper we compared the best known direct formulas and the so-called double loop reordering approach for estimation Sobol' main effect indices on a set of test functions for models with independent and dependent inputs. Both MC and QMC samplings were considered. From the convergence results for models with independent inputs it follows that in majority of tests cases improved direct formulas show much higher efficiency than the original Sobol' formula especially for cases of small values of Sobol' indices with Owen and Oracle formulas outperforming other formulas. DLR outperforms direct formulas on average and by a wide margin when the values of Sobol' indices are not very small. For models with dependent inputs DLR is much easier to implement algorithmically than the direct extended Sobol' formula and hence it is much faster to run. However, practically the DLR method is limited to computing Sobol' main effect indices for a single index only. Convergence of all methods is much higher when QMC sampling is used apart from the case of type C function.

**Acknowledgements**

The authors would like to thank B. Delpuech for his help in preparation of this work. The financial support by the EPSRC grant EP/H03126X/1 is gratefully acknowledged.